\documentclass[11pt,a4paper]{article}
\usepackage{fullpage,amsfonts,amssymb,latexsym,bbm,amsmath}
\usepackage{graphicx} 
\usepackage{indentfirst}
\usepackage{pdfsync}
\usepackage{stackrel}
\usepackage{color}
\usepackage{url}
\usepackage{graphicx}
\usepackage[all]{xy}
\usepackage{amscd}
\usepackage {stackrel}
\usepackage{ifpdf}
\ifpdf%
\else
  \usepackage{pstricks}
  
\fi
\usepackage{mathtools}

\addcontentsline{toc}{section}{APPENDIX}

\newcommand{\RR}{\mathbb{R}}
\newcommand{\CC}{\mathbb{C}}

\newcommand{\ZZ}{\mathbb{Z}}

\newcommand{\NN}{\mathbb{N}}

\newcommand\cH{\mathcal{H}}

\newcommand\bm{\bf{m}}
\newcommand\be{\bf{e}}
\newcommand\bB{\bf{B}}

\def\Re{\mathop{\rm Re} \nolimits}
\def\Im{\mathop{\rm Im} \nolimits}
\def\beginpf{\noindent {\sl Proof:\,\,}}

\newtheorem{lemma}{Lemma}
\newtheorem{corollary}{Corollary}

\newtheorem{proposition}{Proposition}
\newtheorem{remark}{Remark}
\newcommand{\boite}{\mbox{} \hfill \mbox{\rule{2mm}{2mm}}}

\def\endpf{\boite\medskip}

\def\beq{\begin{equation}}
\def\eeq{\end{equation}}

 \title{Solutions to inverse moment estimation problems
      in dimension 2, 
   using best constrained approximation} 
\author{Juliette Leblond\thanks{INRIA, Team Factas, 2004 route des Lucioles, 06902 Sophia Antipolis, France} \and Elodie Pozzi\thanks{Saint Louis University, Department of Mathematics and Statistics, Ritter Hall, 1N. Grand Blvd., Saint Louis, MO 63103, USA}}
\date{}

\begin{document}
\maketitle


\abstract{\noindent We study an inverse problem that consists in estimating the first (zero-order)  moment of 
  some $\RR^2$-valued distribution $\bm$ supported within a closed interval $\bar{S} \subset \RR$,
  from  partial knowledge
  of the solution to the Poisson-Laplace partial differential equation with source term equal to the divergence of $\bm$ 
  on another interval parallel to and located at some distance 
  from $S$. Such a question coincides with a 2D version of an inverse magnetic ``net'' moment recovery question that arises in paleomagnetism,  for thin rock samples. We formulate and constructively solve a best approximation problem under constraint in $L^2$ and in Sobolev spaces involving the restriction of the Poisson extension of the divergence of $\bm$. Numerical results obtained from the described algorithms for the net moment approximation are also furnished.}\\

\noindent {\bf Keywords:} Poisson-Laplace partial differential equation, inverse problems, regularization, best constrained approximation and bounded extremal problems, Poisson and Hilbert transforms, Hardy spaces, Maxwell's equations  and issues.

\section{Introduction}

\noindent 
The present study concerns the practical issue of estimating the ``net'' moment (or the mean value) of some $\RR^2$-valued distribution (or integrable, or square-integrable function) supported on an interval of the real line, 
from the partial knowledge of the divergence of its Poisson extension to the upper half-plane, on some interval located on a parallel line at distance $h>0$ from the first one.

\noindent This is a two dimensional (2D) formulation of a three dimensional (3D) inverse magnetic moment recovery problem in paleomagnetism, for thin rock samples. There, the magnetization is assumed to be  some $\RR^3$-valued distribution (or integrable function) supported on a planar sample (square), of which we aim at estimating the net moment, given measurements of the normal component of the generated magnetic field on another planar measurements set taken to be a square parallel to the sample and located at some distance $h$ to it. Moment and more general magnetization recovery issues are considered in \cite{Betal-IP, BEP-moments, Proc-JE}. 
More about the data acquisition process and the use of scanning SQUID (Superconducting Quantum Interference Device) microscopy devices for measuring a component of the magnetic field produced  by weakly magnetized pieces of rocks can be found in the introductory sections of these references.

\noindent In both situations, the partial differential equation (PDE) model that drives the behaviour of the magnetic field derives  from  Maxwell equations in magnetostatics, see \cite{Jackson}. They ensure that  the magnetic field $\bB$ derives from a  scalar magnetic potential $U$ which is solution to
a  Poisson-Laplace elliptic PDE with right-hand side of divergence form, that relates the Laplacian of that potential to the divergence of the magnetization distribution $\bm$: 
\[
  {\bB} = - \mu_0 \, 
\mbox{grad} \, U \, , \ \Delta \, U = \mbox{div} \, {\bm} \, \text{ in } \RR^n \, ,
  \]
  with $n = 2, 3$, where both $\bB$, $\bm$ are $\RR^n$-valued quantities, and $\mu_0$ is the permeability of the free space. We stick to situations where  $\bm$ has a  compact  support contained in the hyperplane $\{x_n = 0\}$  and measurements of the vertical component $b_n$ of the field ${\bB}$ are available on another compact subset of $\{x_n = h\}$ for some $h >0$. 

\noindent The determination of the magnetic moment (mean value of the magnetization) provides useful preliminary information for the full inversion, i.e. for magnetization estimation, in particular in unidirectional situations where its components are proportional one to the others (the unknown quantity being a $\RR$-valued distribution, the direction/orientation vector being fixed).\\

\noindent In  3D situation,  these issues were efficiently analyzed with tools from harmonic analysis,  specifically Poisson kernel and Riesz transforms \cite{SteinWeiss},  in \cite{Betal-IP, BEP-moments, Proc-JE}, see also references therein, 
using their links with Hardy spaces of harmonic gradients. 
The existence of silent sources in 3D, elements of the non-zero kernel of the non injective magnetization-to-field operator (the operator that maps the magnetization to the measured component of the magnetic field)  was established together with their characterization in \cite{Betal-IP}. This makes non unique the solution to the inverse magnetization issue from the corresponding field data. The mean value of the magnetization however is uniquely determined by the field data.

\noindent In the present 2D case, we will use similar tools, Poisson kernel, Hilbert transform and harmonic conjugation, links with Hardy spaces of functions of the complex variable \cite{Duren, Garnett, Javad};  see also \cite{LebPon, DP} for related issues. 
Considering  square-summable magnetizations $\bm$ supported on an interval $S$ of the real line, 
we will establish that the magnetization-to-field operator ${\bm} \mapsto b_2[{\bm}]$ is injective, whence there are no silent sources, and that the mean value $\langle {\bm} \rangle$ of the magnetization is yet uniquely determined by the field data $b_2[{\bm}]$, provided by values of the vertical 
component on the measurement interval $K$.\\

\noindent Our purpose is to establish the existence and to build linear estimators for the mean values $\langle {m}_i \rangle$ on $S$ of the components $m_i$ of ${\bm}$ for $i=1,2$, as was done in \cite{BEP-moments} for the 3D case. Indeed, though an academic version of the related physical issue, the present 2D situation possesses its own mathematical interest and specificity.
These linear estimators consist in square summable functions $\phi_i$
such that their effect (scalar product) against the data on $K$ is as close as possible to the moment components $\langle {m}_i \rangle$, for all ${\bm}$ bounded by some $L^2$ norm (the functions being built once and for all).
With ${\be}_1=(\chi_S,0)$, ${\be}_2 = (0,\chi_S)$, observe that:
$\langle m_i \rangle = \langle{\bm} \, , \, {\be}_i \rangle_{L^2(S,\RR^2)}$,
whence, if $b_2^*$ denotes the adjoint operator to $b_2$, we have:
\begin{eqnarray*}
  \left\vert \langle b_2[{\bm}] \, , \, \phi_{i} \rangle_{L^2(K)} - \langle m_i \rangle 
  \right\vert
&\leq&\Vert {\bm}\Vert_{L^2(S,\RR^2)}\Vert b_2^*[\phi_{i}]-{\be}_i\Vert_{L^2(S,\RR^2)},
\end{eqnarray*}
\noindent where $\phi_i$is a square integrable function on $K$. We will see that the above minimization problem is still ill-posed, even if uniqueness is granted. Specifically, there exist a sequence of functions such that their scalar product on $K$ against $b_2[{\bm}]$ converges to $\langle m_i \rangle$ as $n \to \infty$. Their quadratic norm however diverge, which reflects an 
unstable behaviour, as is classical in such inverse problems, see \cite{R-al}. Regularization is thus needed (of Tykhonov type), see \cite{Kress} in order to set up and to solve a well-posed problem.

\noindent In order to construct such a numerical magnetometer, we then face the best constrained approximation issue (bounded extremal problem, BEP) of finding the function $\phi_i$ on $K$ satisfying some norm constraint there, such that its scalar product against $b_2[{\bm}]$ is as close as possible to $\langle {m}_i \rangle$.

\noindent Such a norm constraint will be considered both in the Lebesgue space $L^2(K)$ and in the Sobolev space 
$W_0^{1,2}(K) \subset W^{1,2}(K)$.
In those Hilbert spaces, the above problems will be shown to be well-posed, the approximation subsets being closed and convex. The results obtained in $L^2(K)$ and in 
$W_0^{1,2}(K)$ are different. In particular, the solutions in $W_0^{1,2}(K)$ demonstrate less oscillations at the endpoints of $K$ 
than the ones in $L^2(K)$, 
and allow to incorporate  vanishing boundary conditions. The results for moment estimation will be compared between them.

\noindent Preliminary numerical computations in the 3D setting, with planar squared support and measurement set, were run on appropriate finite elements bases, see  \cite{BEP-moments}; they are actually heavy. 
In the present 2D setting and on intervals, they are of course much less costly and we perform some of them using expansions on the Fourier bases, in Lebesgue and Sobolev spaces.\\

\noindent The overview of the present work is as follows. In Section \ref{sec:prel}, we introduce some notation, recall definitions and properties, concerning Lebesgue, Sobolev and Hardy spaces of functions, together with integral Poisson and Hilbert transforms.
We then study in Section \ref{sec:mo} the properties of the operator $b_2 :  {\bm} \to b_2[{\bm}]$, which maps the magnetization to the second (vertical) component of the produced magnetic field, and of its adjoint.
Section \ref{sec:beps} is devoted to the bounded extremal problems of which we consider two versions, in $L^2(K)$ and in $W_0^{1,2}(K)$, establishing their well-posedness and characterizing their solutions. 
Computational algorithms together with results of preliminary numerical simulations are provided in Section \ref{sec:comput}, with figures in an appendix. We finally discuss in Section \ref{sec:pc} some concluding remarks and perspectives.

\section{Preliminaries}
\label{sec:prel}

\noindent Let $\CC_+ \simeq \RR^2_+$ be the upper half-plane $\mathbb{C}_+=\left\{x+iy:y>0\right\} \subset \CC$,  
$\RR^2_+=\left\{(x,y):y>0\right\} \subset \RR^2$. The partial derivatives with respect to $x$ and $y$ will be denoted by $\partial_x$ and $\partial_y$ respectively.
  
\subsection{Lebesgue and Sobolev spaces}
\label{subs:lsp}

\noindent Let $I \subseteq \RR$ a non empty open interval. For $1\leq p<\infty$, we will denote by $L^p(I,\mathbb{K})$ with $\mathbb{K}=\RR$, $\CC$ or $\RR^2$  the space of $\mathbb{K}$-valued functions $f$ such that $\vert f\vert^p$ is integrable on $I$; we will simply write $L^p(I)=L^p(I,\RR)$.

\noindent 
We equip the space $L^2(I, \mathbb{R}^2)$ with the following inner product: for $f=(f_1,f_2)$ and $g=(g_1,g_2)\in L^2(I, \mathbb{R}^2)$, 
\[
\langle f,g\rangle_{L^2(I, \mathbb{R}^2)} =\langle (f_1,f_2),(g_1,g_2)\rangle_{L^2(I, \mathbb{R}^2)}
= \langle f_1,g_1\rangle_{L^2(I)}+\langle f_2,g_2\rangle_{L^2(I)} \, , 
\]
while $L^2(I, \CC)$ is equipped with: 
\begin{equation}
  \label{psIC}
  \langle f,g \rangle_{L^2(I, \CC)}=\langle \Re \, f ,\Re \, g \rangle_{L^2(I)}+\langle \Im \, f ,\Im\,g \rangle_{L^2(I)} \text{ for } f\, , \ g \in L^2(I, \CC) \, .
\end{equation}
For $I \subsetneq \RR$, for $\phi \in L^2(I)$, we will write $\widetilde{\phi}=\phi \vee 0 \in L^2(\mathbb{R})$, the function equal to $\phi$ on $I$ and equal to $0$ on $\RR\backslash I$. Note that we identify $\widetilde{L^2(I)}$ to $\chi_I \,  L^2(\RR)$ and to the subspace of $L^2(\mathbb{R})$ of functions supported on $I$. We also have $\widetilde{\phi}_{|_I} = R_I \, \widetilde{\phi}= \phi$ for $\phi\in L^2(I)$ and the adjoint of the restriction operator $R_I$ from $L^2(\RR)$ to $L^2(I)$ is the extension operator $\widetilde{.}$ on $L^2(I)$.

\noindent When $I$ is bounded, $\widetilde{L^2(I)}\subset L^2(\mathbb{R})\cap L^1(\mathbb{R})$. Moreover, functions in $\widetilde{L^2(I)}$ are compactly supported on $\mathbb{R}$.

\noindent For $I \subseteq \RR$,
 the Sobolev space $W^{1,2}(I)$ is the space of real-valued functions $f \in L^2(I)$ such that their (distributional) derivative $f'$ belongs to $L^2(I)$; $W^{1,2}(I)$ is equipped with the norm (see \cite[Sec. 8.2]{Brezis}):
$$\Vert f\Vert_{W^{1,2}(I)}=\Vert f\Vert_{L^2(I)}+\Vert f'\Vert_{L^2(I)}\, ,\quad f\in W^{1,2}(I)\, .$$
\noindent A function in $W^{1,2}(I)$ can be extended 
to a function in $W^{1,2}(\mathbb{R})$,
see the extension operator in \cite[Thm 8.6]{Brezis}. The space $C(\bar{I})$ denotes the space of continuous functions on $I$ with the norm
$$\Vert f\Vert_{C(\bar{I})}=\sup_{x\in \bar{I}}\vert f(x)\vert \, , \quad  f \in C(\bar{I})\, . $$
\noindent Note that $W^{1,2}(I)\subset C(\bar{I})$, the space of continuous real-valued functions on $\bar{I}$, with compact (hence continuous) injection, see \cite[Thm 8.5]{Brezis}. 
If $I$ is bounded, $C^1(\bar{I}) \subset W^{1,2}(I)$ will denote the space of functions in $C(\bar{I})$ with (strong) derivative in $C(\bar{I})$. The space is equipped with the norm 
$$\Vert f\Vert_{C^1(\bar{I})}=\sup_{x\in \bar{I}}\vert f(x)\vert+\sup_{x\in \bar{I}}\vert f'(x)\vert \, , \quad  f \in C^1(\bar{I})\, ,$$
\noindent where the (strong) derivative coincides with the distributional derivative. The space $C^\infty(\RR)$ will denote the space of functions such that derivatives at all order are continuous. The subspace $W_0^{1,2}(I)$ of $W^{1,2}(I)$ is the collection of functions $f\in W^{1,2}(I)$ such that $f=0$ at the endpoints of $I$. A consequence is that $\widetilde{f}\in W^{1,2}(\mathbb{R})$ for $f\in W_0^{1,2}(I)$; in other words, $\widetilde{.}$ is the extension operator from $W_0^{1,2}(I)$ to $W^{1,2}(\RR)$ (see \cite[Sec 8.3, Remark 16]{Brezis}.
We equip $W_0^{1,2}(K)$ with the norm:
\[\Vert \phi\Vert_{W_0^{1,2}(K)}=\Vert \phi'\Vert_{L^2(K)} \, , \quad  \phi \in {W_0^{1,2}(K)}\, ,\]
which is 
equivalent to $\Vert  \phi \Vert_{W^{1,2}(K)}$
by Poincar\' e's inequality (see \cite[Prop 8.13]{Brezis}).

\subsection{Poisson kernel, conjugate Poisson kernel and Hilbert transform}
\label{subs:pch}
\noindent For $y>0$, we denote by $P_y$ the Poisson kernel of the upper half-plane and by $Q_y$ the conjugate Poisson kernel. For $x \in \RR$:
  $$ P_y(x)=\frac{1}{\pi}\frac{y}{x^2+y^2},\quad Q_y(x)=\frac{1}{\pi}\frac{x}{x^2+y^2}. $$
  \noindent see \cite{Duren, Garnett, SteinWeiss}. Let $x \in \RR$, $y>0$. Then
$$ P_y\star u(x)=\frac{1}{\pi}\int_{-\infty}^{\infty}\frac{y}{(x-t)^2+y^2}u(t) dt,\,\,\, u\in L^2(\mathbb{R}),$$
  and $\mathcal{P}_y:u \mapsto P_y\star u$ is a bounded (real linear) operator from $L^2(\mathbb{R})$   to $L^2(\mathbb{R})$ (\cite[Thm 1.3]{SteinWeiss}).

\noindent The map, $x \mapsto P_y(x)$  belongs to $L^1(\RR)$ (where its norm is equal to 1). Moreover,  $(x,y) \mapsto P_y(x)$  belongs to $C^\infty(\RR_+^2)$. Therefore, its partial derivatives  are also $C^\infty(\RR_+^2)$-smooth functions, whence for $y>0$,  $x \mapsto (\partial_x P_y)(x)$ and $x \mapsto (\partial_y P_y)(x)$ belong to $C^\infty(\RR)$. The functions $P_y$ and $Q_y$ satisfy the Cauchy-Riemann equations on $\CC_+$: $\partial_yP_y=-\partial_xQ_y$, $\partial_xP_y=\partial_yQ_y$.\\

\noindent It can be checked by direct computations that $x \mapsto (\partial_x P_y)(x)$ and $x \mapsto (\partial_y P_y)(x)$ belong to $L^1(\RR)$, and $\|\partial_x P_y\|_{L^1(\RR)}=\|\partial_y P_y\|_{L^1(\RR)}=2 / (\pi \, y)$. Hence, $u \to (\partial_x P_y)\star u$ and $u \to (\partial_y P_y)\star u$ are  bounded operators from $L^2(\mathbb{R})$   to $L^2(\mathbb{R})$. Further, for $u \in L^2(\mathbb{R})$, $(x,y) \mapsto \mathcal{P}_y [u](x)=(P_y\star u) (x)$ is also harmonic  in $\RR_+^2$, hence a $C^\infty(\RR_+^2)$-smooth function. Thus, its partial derivatives $x \mapsto \partial_x( P_y \star u )(x)$ and $x \mapsto \partial_y( P_y \star u )(x)$ belong to $C^\infty(\RR)$, for $y>0$.\\
 
\noindent Using the Fourier transform as in \cite[Sec. 2.2]{Proc-JE}, we see that convolution by the Poisson kernel and differentiation commute: for $y>0$ and $u \in L^2(\mathbb{R})$, we have that $\partial_x (P_y\star u) = (\partial_x P_y) \star u$.\\ 
In particular, the operator $\mathcal{P}_y $ maps continuously $L^2(\mathbb{R})$ onto $W^{1,2}(\mathbb{R})$. By Cauchy-Riemann equations, $Q_y$ satisfies the  properties of $P_y$ described above. Indeed, $u\mapsto \partial_x(Q_y\star u)$ and $u\mapsto \partial_y(Q_y\star u)$ are bounded on $L^2(\mathbb{R})$ and for $y>0$, $x\mapsto \partial_x(Q_y\star u)$ and $x\mapsto \partial_y(Q_y\star u)$ are $C^\infty(\RR)$.\\

\noindent In order to simplify the notations, we may not use any parentheses in expressions like $\partial_x P_y \star u$, $\partial_y P_y \star u$, $\partial_x Q_y \star u$ and $\partial_y Q_y \star u$ for $u \in L^2(\RR)$ and $y>0$.\\

\noindent Let $I, \, J \subsetneq \RR$ be non empty open bounded interval. 
Then, for $y >0$, $R_I \, {\cal P}_y$ is a Hilbert-Schmidt operator on $\widetilde{L^2(J,\CC)} \subset L^2(\RR,\CC)$ (see \cite[p. 264, Ex. 2.19]{Kato}). 
Indeed, for $u\in L^2(J)$, we have
$$R_I \, {\cal P}_y[u](x) = R_I \, \left(P_y\star u\right)(x)=\frac{1}{\pi}\int_J \frac{h}{(x-t)^2+y^2}u(t) dt,\,\,x\in I \, ,$$
and $R_I\, {\cal P}_y$ is an integral operator on $\widetilde{L^2(J,\CC)}$ with kernel $k(x,t)=P_y(x-t)$, $k\in L^2(I \times J)$, since ${P}_y(x)$ is a smooth function for $y>0$. 
Thus, it is compact. The same conclusion holds for the operator $u \mapsto R_I \,  \left(\partial_x \, P_y \star u\right)$ (and $u \mapsto R_I \,  \left(\partial_y \, P_y\star u\right)$) on $\widetilde{L^2(J,\CC)}$.\\

\noindent For $u\in L^2(\mathbb{R})$, the functions $\mathcal{P}_y [u]=P_y\star u$ and $Q_y\star u$ have non-tangential limits as $y$ tends to zero and 
\begin{equation}\label{ntlimit}
\lim_{y\to 0}P_y\star u(x)=u(x)\,\,\,\hbox{and}\,\,\,\lim_{y\to 0}Q_y\star u(x)=\mathcal{H}u(x),
\end{equation}
\noindent where $\mathcal{H}$ is the Hilbert transform from $L^2(\mathbb{R})$ to $L^2(\mathbb{R})$ defined by 
$$\mathcal{H}u(x)=\displaystyle{\frac{1}{\pi}\lim_{\varepsilon\to 0}\int_{\vert x-t\vert>\varepsilon}\frac{u(t)}{x-t}dt}.$$
\noindent The Hilbert transform is bounded and isometric on $L^2(\mathbb{R})$ and satisfies $\mathcal{H}^2=-I$, whence its adjoint $\mathcal{H}^* = - \mathcal{H}$. The Hilbert transform commutes with the Poisson operator (which can be proved using the Fourier transform): 
\begin{equation}\label{com:hilbert-poisson}
 \cH (P_y\star u)=\cH (\mathcal{P}_y[u])=P_y\star \cH u =\mathcal{P}_y[\cH u ] \, .
 \end{equation}
\noindent  For $y>0$, $x\in\mathbb{R}$ and $u\in L^2(\mathbb{R})$, one can prove using \cite[Thm 11.6]{Javad} that:
\begin{equation}\label{eg1}
Q_y\star u ={P}_y\star \mathcal{H} u=\mathcal{P}_y[\cH u ]\, .
\end{equation}
\noindent Hence, from \eqref{com:hilbert-poisson} and \eqref{eg1}, we get that $Q_y =\mathcal{H} \, {P}_y$ on $L^2(\RR)$. We will also use in the sequel the following equality, for $g \in L^2(\RR)$:
\begin{equation}\label{poisson_der}
\partial_y \mathcal{P}_y[g]  = -\mathcal{H} \left(\partial_x \mathcal{P}_y [g] \right).
\end{equation}

\noindent It can be obtained using the Fourier transform 
as in \cite{Proc-JE, SteinWeiss}. Actually,  for $g\in L^2(\RR)$, $\partial_y \, \mathcal{P}_y[g] =-\left[\partial_x \, Q_y \star g\right]$, and $\partial_x \, \mathcal{P}_y[g]= \partial_y \,Q_y \star g$. And as a consequence of \eqref{poisson_der}, and the isometric character of $\mathcal{H}$, we get $\|\partial_y \mathcal{P}_y[g]\|_{L^2(\RR)}= \|\partial_x \mathcal{P}_y[g]\|_{L^2(\RR)} $.\\

\noindent We also have the following uniqueness result, from \cite[Lem. 2.19]{Proc-JE}:
\begin{lemma}
  \label{lem:ph}
  Let $h>0$ and $g \in L^2(\RR)$ such that $\partial_x \, \mathcal{P}_h[g] = 0$ on a non-empty open subset 
  of $\RR$. Then $g \equiv 0$ on $\RR$.
\end{lemma}
%
As a corollary to Lemma \ref{lem:ph}, we get that if $g\in L^2(\RR)$ is such that $\mathcal{P}_h [g] = P_h\star g=0$ on a non-empty open subset of $\RR$ for some $h>0$, then $g\equiv 0$ on $\RR$.
%

\subsection{Hardy space $H^2$}\label{sec:hardy}
\label{subs:hs}
      
\noindent For $u\in L^2(\RR)$, the complex-valued function $P_y\star u+iQ_y\star u$ belongs to $H^2(\CC_+)=H^2(\CC_+, \CC)$, the space of analytic functions $F$ on $\CC_+$ such that
\begin{equation}\label{hardy:ana}
\Vert F\Vert_{H^2(\mathbb{C}_+)}:=\sup_{y>0}\left(\int_{-\infty}^{\infty} \vert F(x+iy)\vert^2 dx \right)^{1/2}<\infty \, ,
\end{equation}
\noindent see \cite{Duren, Garnett}. A function $F\in H^2(\CC_+)$ admits a non-tangential limit $F^*\in L^2(\RR, \CC)$ as $y$ tends to $0$, almost everywhere on $\RR$. 
The function $F^{*}$ belongs to the space $H^2(\RR)=H^2(\RR, \CC)$ of 
functions $f \in L^2(\RR)$ such that 
$$\int_{-\infty}^{\infty}\frac{f(t)}{t-\overline{z}}dt=0,\,\,\,z\in\CC_+. $$
\noindent The space $H^2(\RR)$  coincides with the space of boundary values  on $\RR$ (also ``traces'', if smoothness allows) of $H^2(\CC_+)$ functions.
Indeed, the map $f\in H^2(\RR)\longmapsto P_y\star f$ is an isometric isomorphism from $H^2(\RR)$ onto $H^2(\CC_+)$.\\
It is not difficult to see that $(I+i\mathcal{H})\,L^2(\RR)=H^2(\RR)$. As a consequence, the map $u\longmapsto P_y\star(I+i\mathcal{H})u$ is an isomorphism from $L^2(\RR)$ onto $H^2(\CC_+)$. A function $f\in H^2(\CC_+)$ satisfies the boundary uniqueness Theorem \cite[Cor. 6.4.2]{Niko}: if $F\in H^2(\RR)$ is such that $F=0$ on $I\subset 
\RR$ with $|I|>0$ then $F\equiv 0$ on $\RR$.

\noindent Note that if  $I\subset\RR$, the space $(I+i\mathcal{H})(\widetilde{L^2(I)})$ coincides with the space of functions $F^{*}$ of $H^2(\RR)$ such that $\text{supp} \Re \, F^{*} \subset I$.
         
\noindent Observe that for $f\in L^2(\RR,\CC)$,
$$(I+i\cH) \, f =2\,(I+i\cH)\, [\Re \, P_+(f)] =2 \, P_+(f),$$
\noindent where $P_+(f)$ is the orthogonal projection from $L^2(\RR,\CC)$ onto $H^2(\RR)$, see \cite{LebPon}.\\

\noindent Finally, an amazing property of the Hilbert transform is given by the  following Lemma (of which a proof can be found in \cite[Lem. 2.12]{R-al}): 
\begin{lemma}\label{lem:hilbtronc}
If $I$, $J\subset\RR$ with $|I|,\,|J|>0$ such that $I \cap J =\emptyset$, then for any $g\in L^2(\RR)$ such that $\cH(\chi_Jg)=0$ on $I$, we have $g=0$ on $J$. 
\end{lemma}
Lemma \ref{lem:hilbtronc} can be seen as a consequence of the boundary uniqueness Theorem in $H^2(\RR)$. 
Indeed, if $\cH(\chi_J g)=0$ on $I$ then the function in $H^2(\RR)$ given by $F=\chi_J g+i\cH(\chi_J g)$ is equal to $0$ on $I$ (as $J\cap I=\emptyset$). By the boundary uniqueness Theorem, it follows that $F=0$ on $\RR$ and $g=0$ on $J$.

\noindent When $I\cap J\neq \emptyset$ and $I \not \subseteq J$, Lemma \ref{lem:hilbtronc} still holds true in some cases:
\begin{enumerate}
\item[$\bullet$] $J\subsetneq I$: indeed, the function $F=\chi_J g+i\cH(\chi_J g)$ is equal to zero on $I\backslash J$ and the boundary uniqueness Theorem implies that $g=0$ on $J$, whenever $|I\backslash J| >0$;
\item[$\bullet$] $I\cap J\neq\emptyset$ and $I\not\subset J$: the same function $F$ is equal to zero on $I\backslash (I\cap J)$ and $g=0$ on $J$ by the boundary uniqueness Theorem, provided that $|I\backslash (I\cap J)|>0$. 
\end{enumerate}
\noindent However, if $I \subseteq J$ are bounded intervals, then there exists a function $g\in L^2(\RR)$, $g\neq 0$, such that $\cH(\chi_Jg)\neq 0$. 
With $I=(a,b)$, the (non-zero) function defined on $\RR$ by $g(x)=\chi_I(x)/{\sqrt{(x-a)(b-x)}}$ is such that $\cH(\chi_I g)=\cH(g)=0$ on $I$. Indeed, the function $x \mapsto {1}/{\sqrt{(x-a)(b-x)}}$ is shown in \cite[Lem. 2.1]{kop} to be an eigenfunction of the operator $\cH(\chi_I \cdot)_{|_I}$ associated to the eigenvalue 0.
\begin{enumerate}
\item[$\bullet$] \noindent If  $I=J$, the function $g$ furnishes a solution to the issue.
\item[$\bullet$] \noindent If $I\subsetneq J$, one can write 
$\cH(\chi_Jg)=\cH(\chi_Ig)+\cH(\chi_{J\backslash I} g)$.
\noindent 
The above function $g$ vanishes outside $I$ hence on $J\backslash I$.
Thus, $\cH(\chi_{J\backslash I} g)=0$ on $\RR$. 
Further, $\cH(\chi_Ig)=0$ on $I$, whence $\cH(\chi_J g)=0$ on $I$, 
but $g\neq 0$ on $J$.
\end{enumerate}
\section{Main operators}
\label{sec:mo}
\noindent Let $S$ and $ K \subset \RR$ be two nonempty open bounded intervals.   Fix $h > 0$. Let ${\bm} \in L^2(S, \mathbb{R}^2)$. 

\noindent Taking the convolution of  the PDE $\Delta \, U = \mbox{div} \,  {\bm}$ by the fundamental solution to Laplace equation in dimension $n=2$ and applying Green formula, we obtain, for $X \in \mathbb{R}^2 \setminus \text{supp} \, {\bm}$:
\[
U[{\bm}](X) = -\frac{1}{2 \, \pi} \, \iint_{\RR^2} \left(\mbox{div}_Y \, {\bm}\right) (Y)\, \log |X-Y| \, d \, Y
= \frac{1}{2 \, \pi} \, \iint_{\RR^2} {\bm} (Y) \, \cdot \, \mbox{grad}_Y \log |X-Y| \, d \, Y 
\]
\[
= \frac{1}{2 \, \pi} \, \iint_{\RR^2} \frac{{\bm} (Y)\cdot \left(X-Y\right)}{|X-Y|^2} \, d \, Y \, . 
\]
Since the support of ${\bm}$ is a subset of $S \times \{0\}$, we get, at $X = (x,y)\in \mathbb{R}^2 \setminus S \times \{0\}$, with $Y = (t,0)$:
\[
U[{\bm}](x,y) = \frac{1}{2 \, \pi} \, \int_{S} \frac{m_1(t) (x-t) + m_2(t)y}{(x-t)^2 + y^2} \, d \, t = \frac{1}{2} \, \left(Q_y \star\widetilde{m_1}+{P}_y \star \widetilde{m_2}\right)(x,y)\, .
\]
\noindent Let then  $b_2[{\bm}] = -\left(\partial_y\,U[{\bm}]\right)_{|_{K\times\{h\}}}$ (from now, we will ignore  on the multiplicative constant ${\mu_0}/{2}$): 
\begin{eqnarray*}
b_2: & L^2(S, \mathbb{R}^2) & \longrightarrow \ L^2(K)\\
& {\bm}=(m_1,m_2)&\longmapsto \ -\partial_y\,\left(Q_y \star\widetilde{m_1}+{P}_y \star \widetilde{m_2}\right)_{|_{K\times\{h\}}}\, ,
\end{eqnarray*}
with the use of Cauchy-Riemann equations, which also allow to get 
\begin{eqnarray*} 
  b_2[\bm]&=&  -\left(\partial_x\, \mathcal{P}_y [\widetilde{m_1}] +\partial_y\,\mathcal{P}_y[ \widetilde{m_2}]\right)_{|_{K\times\{h\}}}= -\left[\partial_x\, \left(P_y \star \widetilde{m_1} -Q_y\star \widetilde{m_2}\right)\right]_{|_{K\times\{h\}}} \\
&=&-\partial_x \mathcal{P}_y[\widetilde{m_1}-\mathcal{H}\widetilde{m_2}]_{|_{K\times\{h\}}}\, ,
\end{eqnarray*}
\noindent where the second equality follows from \eqref{eg1}. Using the boundedness of the convolution by ${P}_y$ from $L^{2}(\mathbb{R})$ to $W^{1,2}(\mathbb{R})$, the boundedness of the Hilbert transform and the equality $\Vert\widetilde{m_i}\Vert_{L^2(\mathbb{R})}=\Vert m_i\Vert_{L^2(S)}$, for $i\in\{1,2\}$ (see Section \ref{sec:prel}), one can prove that the operator $b_2$ is bounded $L^2(S,\mathbb{R}^2)\rightarrow L^2(Q)$.  

\begin{proposition}
  \label{propb2star}
The adjoint operator $b_2^*$  to $b_2$ is given by 
\begin{eqnarray*} 
b_2^{*}: & L^2(K) & \longrightarrow  L^2(S,\mathbb{R}^2)\\
& \phi &\longmapsto (\partial_x \, (P_y\star \widetilde{\phi})_{|_{S\times\{h\}}},\partial_x \, (Q_y\star\widetilde{\phi})_{|_{S\times\{h\}}})
\end{eqnarray*}
and is such that
\[
b_2^{*}[\phi]=(((\partial_x \, P_y)\star \widetilde{\phi})_{|_{S\times\{h\}}},((\partial_x \, Q_y)\star\widetilde{\phi})_{|_{S\times\{h\}}}) 
=\partial_x \, (P_y\star(I,\mathcal{H})\widetilde{\phi})_{|_{S\times\{h\}}}.
\]
\end{proposition}

\beginpf \noindent Let $y>0$ fixed. For ${\bm}\in L^2(S, \mathbb{R}^2)$ and $\phi\in L^2(K)$, we have that 
\begin{eqnarray*}
 & & \langle b_2[{\bm}],\phi\rangle_{L^2(K)} = -\langle\partial_x(P_y\star\widetilde{m_1}),\phi \rangle_{L^2(K)}+\langle\partial_x(Q_y\star\widetilde{m_2}),\phi \rangle_{L^2(K)}\\
 & &  =-\int_\RR \int_\RR \partial_x(P_y(x-t))\widetilde{m_1}(t)\widetilde{\phi}(x) dtdx+\int_\RR\int_\RR \partial_x(Q_y(x-t))\widetilde{m_1}(t)\widetilde{\phi}(x) dtdx\\
 & & =\int_\RR \left(\int_\RR \partial_t(P_y(x-t))\widetilde{\phi}(x) dx\right) \widetilde{m_1}(t)dt+\int_\RR\left(\int_\RR \partial_t(Q_y(x-t))\widetilde{\phi}(x)dx\right)\widetilde{m_1}(t)dt\\
 & & =\int_S \partial_t\left(\int_\RR P_y(t-x)\widetilde{\phi}(x) dx\right) m_1(t)dt+\int_S\partial_t\left(\int_\RR Q_y(t-x) \widetilde{\phi}(x) dx\right) m_1(t)dt\\
 & &  =\langle m_1,\partial_t(P_y\star \widetilde{\phi})\rangle_{L^2(S)}+\langle m_2,\partial_t(Q_y\star \widetilde{\phi})\rangle_{L^2(S)}\\
 & &  =\langle (m_1,m_2),(\partial_t(P_y\star \widetilde{\phi}),\partial_t(Q_y\star \widetilde{\phi}))\rangle_{L^2(S,\RR^2)}, \, 
\end{eqnarray*}

\noindent with 
\endpf\\
These operators have null kernel, as follows from the next results.
\begin{proposition}
  \label{lem:b2inj}
The operator  $b_2$ is injective. 
\end{proposition}
\beginpf
Let ${\bm} \in \text{Ker} \, b_2$. 
Then, $\partial_x P_y \star(\widetilde{m_1}-\mathcal{H}\widetilde{m_2}) = 0$ on $K\times \{h\}$, which implies that  $\widetilde{m_1}-\mathcal{H}\widetilde{m_2} = 0$ on $\RR$ in view of Lemma \ref{lem:ph}. 
Hence, $\widetilde{m_1}=\mathcal{H}\widetilde{m_2}$ which implies that $\mathcal{H}\widetilde{m_2} = 0$ outside $S$ whence so does $\widetilde{m_2}+i\mathcal{H}\widetilde{m_2}$. Because it coincides with the boundary value of a function in $H^2(\mathbb{C}_+)$, namely $P_y \star (\widetilde{m_2}+i\mathcal{H}\widetilde{m_2})$, this implies that  $m_2=0$ (see \cite{Duren, Garnett}). Therefore, $m_1=0$ and $\text{Ker} \, b_2 = \{0\}$. \endpf\\
\noindent In particular, contrarily to the situation in dimension 3 (and in higher dimensions) see \cite{Betal-IP}, there are no non-vanishing ``silent'' sources (the corresponding so-called ``null-space'' is reduced to $\{0\}$).
\begin{proposition}\label{prop:b2*inj}
The operator $b_2^*$ is injective. 
\end{proposition}
\beginpf Let $f \in \text{Ker} \, b_2^*$. Then, we have that $\partial_x(P_y\star (I,\cH)\widetilde{f})=(0,0)$ on $S\times\{h\}$ which implies that $\partial_x(P_h\star\widetilde{f})=\partial_x(P_h\star \cH\widetilde{f})=0$ on $S$. By Lemma \ref{lem:ph}, it follows that $\widetilde{f}=0$ on $\RR$ and $f=0$ on $S$.
 \endpf\\
 Moreover:
 \begin{corollary}
   \label{cor:denserange}
The operators $b_2$ and $b_2^*$ have a dense range. 
\end{corollary}
\beginpf Since $b_2^*$ is continuous (because so is $b_2$), the following orthogonal decomposition holds true:
\[
L^2(S, \RR^2) =  \text{Ker} \, b_2 \oplus \overline{\text{Ran} \, b_2^*} \, .
\]
Hence, using Proposition \ref{prop:b2*inj}, $b_2^*$ has dense range in $L^2(S, \RR^2)$.
Similarly, 
\[
L^2(K) =  \text{Ker} \, b_2^* \oplus \overline{\text{Ran} \, b_2} \, .
\]
Hence, using Proposition \ref{lem:b2inj}, $b_2$ has dense range in $L^2(K)$. 
\endpf
\begin{remark}
\label{rmkB2}
  \noindent Using the following equality for $f \in L^2(\RR, \CC)$:
\begin{equation}\label{dim:eq}
(I+i\cH)\,f=(I+i\cH)\,(\Re\, f -\cH \, \Im \, f).
\end{equation}
with $\mathcal{H} u = \mathcal{H} [- i^2 u] = - i \mathcal{H} [i u] $ whence  $ \mathcal{H} [i u] = i \mathcal{H} u$, we could consider the analytic versions of $b_2$ and $b_2^{*}$ denoted by $B_2$ and $B_2^{*}$ respectively (see \cite{LebPon}):
$$\begin{matrix}
B_2:L^2(S, \CC)&\longrightarrow &L^2(K)\\
\quad m=m_1+i m_2&\longmapsto&-\Re\partial_x(P_y\star(I+i\mathcal{H}) \widetilde{m})_{|_{K\times\{h\}}}
\end{matrix}, $$ 
$$\begin{matrix}
B_2^{*}:L^2(K)&\longrightarrow &L^2(S, \CC)\\
\quad \phi&\longmapsto&\partial_x(P_y\star(I+i\mathcal{H}) \widetilde{\phi})_{|_{S\times\{h\}}}
\end{matrix}.$$ 
\noindent Let $J$ be the isomorphism from $L^2(S, \CC)$ onto $L^2(S, \RR^2)$ defined by $J(\phi_1+i\phi_2)=(\phi_1,\phi_2)$ which is unitary when $L^2(S, \CC)$ is equipped with the inner product \eqref{psIC}.
%
%
We  have that $B_2=b_2\circ J$ and $b_2^{*}[\phi]=J\circ B_2^{*}[\phi]$, $\phi\in L^2(K)$.
%
As a consequence, $B_2$ and $B_2^*$ are injective and $\overline{\hbox{Ran} \,B_2}=L^2(K)$, $\overline{\hbox{Ran} \,B_2^*}=L^2(S,\CC)$. 
\end{remark}
\begin{proposition}\label{prop:ran}
The ranges of $b_2$ (respectively $B_2$) and $b_2^*$ (respectively $B_2^*$) are not closed: ${\hbox{Ran} \,b_2} \subsetneq L^2(K)$, ${\hbox{Ran} \,b_2^*}\subsetneq L^2(S,\RR^2)$. 
\end{proposition}

\beginpf If we assume that the range of $B_2^*$ is closed, then Corollary \ref{cor:denserange} and Remark \ref{rmkB2} are to the effect that $B_2^*$ is surjective. Hence, for $I\subsetneq S$, with $|I|>0$, there exists $\phi\in L^2(K)$ such that 
\begin{equation}\label{eq:notsurj}
\partial_x(P_h\star(I+i\cH)\widetilde{\phi})_{|_S}=\chi_I \, ,
\end{equation}

\noindent It follows that $\partial_x(P_h\star(I+i\cH)\widetilde{\phi})=0$ on $S\backslash I$ whence $(I+i\cH)\widetilde{\phi}=0$ on $\RR$ by Lemma \ref{lem:ph}. As in the proof of Proposition \ref{lem:b2inj}, this leads to $\phi=0$ which contradicts Equation \eqref{eq:notsurj}. We conclude that $B_2^*$ is not surjective. As $b_2^*=J\circ B_2^*$, we get that $b_2^*$ is not surjective.

\noindent Now, if the range of $B_2$ is closed, then for any $g\in L^2(K)$, there would be $\phi\in L^2(S,\CC)$ such that 
\begin{equation}\label{not-surj}
g=-\Re \partial_x(P_h\star (I+i\cH)\widetilde{\phi})_{|_{K\times\{h\}}}.
\end{equation}

\noindent Let $I\subsetneq K$ with $|I|>0$ and $g=\chi_I$. Then, there exists $\phi\in L^2(S,\CC)$ such that 

$$\Re \partial_x(P_h\star (I+i\cH)\widetilde{\phi})_{|_K}=\partial_x\left(P_h\star \Re\,(I+i\cH)\widetilde{\phi}\right)_{|_K}=0\quad\hbox{on } K\backslash I. $$

\noindent By Lemma \ref{lem:ph}, we get that 
$\Re\,(I+i\cH)\widetilde{\phi}=0\quad\hbox{on }\RR$,
and $\Re\widetilde{\phi}-\cH\Im\widetilde{\phi}=0$ on $\RR$. This means that $\cH\Im\widetilde{\phi}=\Re\widetilde{\phi}$ on $\RR$ whence $\cH\Im\widetilde{\phi}=0$ on $\RR\backslash S$, which proves that $\Im\widetilde{\phi}=0$ on $\RR$ by Lemma \ref{lem:hilbtronc}, and thus $\Re\widetilde{\phi}=0$ on $\RR$. We conclude that $\phi=0$ on $\RR$ which contradicts Equation \eqref{not-surj}. It follows that the range of $B_2$ is not closed and thus $b_2=J^{-1}\circ B_2$ is not surjective.  
\endpf
\begin{remark} \label{ranb2sm}
The ranges of  $b_2$ (and  $B_2$) and $b_2^*$ (and $B_2^*$) are actually made of $C^\infty$-smooth functions on $K$ and $S$, respectively, see Section \ref{subs:pch}.
\end{remark}

\begin{proposition}\label{comp}
The operators $b_2$, $B_2$, $b_2^*$ and $B_2^*$ are compact. 
\end{proposition}

\beginpf 
Because $u \mapsto \left(\partial_x \, P_y \star u\right)_{|_S}$ and $u \mapsto \left(\partial_y \, P_y\star u\right)_{|_S}$ are compact operators on $\widetilde{L^2(K,\CC)}$ for $y>0$ (see Section \ref{subs:pch}), so is $B_2^*=\left(\partial_x \, {\cal P}_h \, (I+iH)\right)_{|_S}$,
and thus also $B_2$.
\noindent As $b_2^*=J\circ B_2^*$ and $b_2=B_2\circ J^{-1}$, we get that $b_2$ and $b_2^*$ are compact.  
\endpf
\begin{remark}
  \noindent
We will also make use of the operator  $a_2$  defined by
$$\begin{matrix}
a_2:L^2(S,\mathbb{R}^2)&\longrightarrow & L^{2}(K)\\ 
\quad {\bm}=(m_1,m_2)&\longmapsto&\left[P_y  \star(\widetilde{m_1}-\mathcal{H}\widetilde{m_2}) \right]_{|_{K\times\{h\}}}
\end{matrix} \, ,$$ 
which is such that $b_2[{\bm}]= -\partial_x (a_2[{\bm}])$.

\noindent
We deduce from the proof of Proposition \ref{propb2star} that $a_2^*[\phi]=(P_y\star \widetilde{\phi}, P_y\star \cH\widetilde{\phi})_{|_{S\times\{h\}}}$, for $\phi\in L^2(K)$, and from the one of Proposition \ref{comp}  that $a_2$ and $a_2^*$ are Hilbert-Schmidt operators.  
Further, whenever $\phi \in W_0^{1,2}(K)$, then $b_2^*[\phi] = -a_2^*[\partial_x \phi]$.
\end{remark}

\section{Bounded extremal problems for moments estimates}
\label{sec:beps}

\noindent We are interested in solving the following bounded extremal problem (BEP), with ${\be} \in L^2(S, \RR^2)$ and $M\geq 0$: find  $\phi_{o} \in {L^2(K)}$,  $\Vert \phi_{o} \Vert_{L^2(K)}\leq M$ that satisfies
\[
\Vert b_2^{*}[\phi_{o}]-{\be}\Vert_{L^2(S,\RR^2)}=\min_{\phi\in L^2(K)\atop{
\Vert \phi\Vert_{L^2(K)}\leq M}}\Vert b_2^{*}[\phi]-{\be}\Vert_{L^2(S,\RR^2)}\label{BEP}\tag{BEP$_{mo}$}
\]
\noindent Our motivation is as follows.  Let ${\be}_1=(\chi_S,0)$, ${\be}_2 = (0,\chi_S)$. A solution $\phi_o$ to \eqref{BEP}  furnishes a linear estimator in $L^2(K)$ for net moment estimates: indeed, for $i\in\{1,2\}$,  
\[ \left\vert \langle b_2[{\bm}] \, , \, \phi_{i} \rangle_{L^2(K)} - \langle m_i \rangle 
  \right\vert=\left\vert\langle{\bm} \, , \, b_2^{*}[\phi_{i}] - {\be}_i \rangle_{L^2(S,\RR^2)}\right\vert\\
  \leq \Vert {\bm}\Vert_{L^2(S,\RR^2)}\Vert b_2^*[\phi_{i}]-{\be}_i\Vert_{L^2(S,\RR^2)}\, ,\]
\noindent where $\phi_{i}=\phi_o[{\be}_i]$ is the solution of \eqref{BEP} with ${\be}={\be}_i$ and $\langle m_i \rangle=\int_S m_i(t) dt$.\\

\noindent Below, we  consider BEP both in the Lebesgue space $L^{2}(K)$ and in the Sobolev space $W_0^{1,2}(K)$ in order to control the derivative and the oscillations of the solution.

\subsection{BEP in $L^{2}(K)$}

\subsubsection{Well-posedness} 

\begin{proposition}
  \label{prop:wpBEP}
  Let  ${\be} \in L^2(S, \RR^2)$ and $M\geq 0$. 
 Problem \eqref{BEP} 
 admits a unique solution $\phi_o \in L^2(K)$. Moreover, if
 $ {\be}  \not \in  \text{Ran}\,b_2^*$, then $\phi_o$ saturates the constraint: $\Vert \phi_o \Vert_{L^2(K)}= M$. 
\end{proposition}
Note that $\Vert \phi_o \Vert_{L^2(K)}= M$  also if ${\be}= b_2^{*}[\phi]$ for some  $\phi \in L^2(K)$ such that $\Vert \phi\Vert_{L^2(K)}\geq M$.\\
\mbox{}\\
\beginpf 
Since $\text{Ran} \, b_2^*$ is dense in $L^2(S, \RR^2)$, the assumptions of \cite[Lem. 2.1]{CP} (see also \cite{CPbook} and \cite[Lem. 1]{LebPon}) are satisfied, to the effect that if $M >0$, Problem 
  \eqref{BEP} admits a solution, which is unique and saturates the constraint if ${\be} \not \in \text{Ran} \, b_{2}^*$ (or if ${\be} =b_{2}^*[\phi]$ for $\phi \in {L^2(K)}$ such that $\Vert \phi\Vert_{L^2(K)}\geq M$). 
  Observe further that because $b_2^*$ is injective (from Proposition \ref{prop:b2*inj}), 
  uniqueness still holds when ${\be} \in \text{Ran} \, b_{2}^*$.
 \endpf\\

  \noindent Existence and uniqueness could also be established by projection onto the closed convex set \cite[Thm 5.2]{Brezis} $b_2^*\left[B_M\right]$,
for $  B_M=\left\{\phi \in L^2(K) \, , \ \Vert \phi\Vert_{L^2(K)}\leq M\right\}$. 
Indeed, $B_M$ is weakly compact as a ball in $L^2(K)$ and 
the weak closedness property of $b_2^*\left[B_M\right]$ follows from the continuity of $b_2^*$, see \cite[Prop. 1]{BEP-moments}. By convexity of $B_M$ and continuity of $b_2^*$, $b_2^*[\phi_n]$ strongly converges to $b_2^*[\phi]$ in $L^2(S, \RR^2)$ whenever $\phi_n$ weakly converges to $\phi$ in $L^2(K)$. Note that the compactness of $b_2^*$ is enough to get that $b_2^*\left[B_M\right]$ is closed. That the constraint is satisfied if ${\be} \not \in b_2^*\left[B_M\right]$ could be established by differentiation of the squared criterion.
\begin{remark} \label{rmk:ip1}
  The density property of $b_2^*\left[L^{2}(K)\right]$ in $L^2(S,\RR^2)$ implies that without the norm constraint, Problem \eqref{BEP} above is ill-posed. Indeed, if ${\be}\not\in b_2^*\left[L^{2}(K)\right]$ the density gives the existence of a sequence $(\phi_n)$ of functions in $L^2(K)$ such that $\Vert b_2^{*}[\phi_n]-{\be}\Vert_{L^2(S,\RR^2)} \to 0$. We claim that $\Vert \phi_n \Vert_{L^2(K)} \to \infty$ as $n \to \infty$. This can be seen by assuming $\Vert \phi_n \Vert_{L^2(K)}$ to remain bounded and extracting a weakly convergent subsequence, say $(\phi_n)$ again, to some $\phi \in b_2^*\left[B_M\right]$ (the later being weakly closed as we saw in the proof of Proposition \ref{prop:wpBEP}). By convexity of $b_2^*\left[B_M\right]$, it follows that $(\phi_n)$ converges strongly to $\phi$ in $L^2(S, \RR^2)$ and this implies that ${\be}= b_2^{*}[\phi] \in \text{Ran}\, b_2^{*}$ which leads to a contradiction. 
\end{remark}

\subsubsection{A constructive solution} 
\label{sec:csol}

\noindent Assume that $ {\be}  \not \in  \text{Ran}\,b_2^*$. From Proposition \ref{prop:wpBEP}, the solution of (BEP$_{mo}$) is a minimum on $L^2(K)$ of the functional $\phi \mapsto \Vert b_2^*[\phi]-{\be}\Vert_{L^2(S,\RR^2)}^2$ which satisfies  the constraint $\Vert \phi\Vert^2_{L^2(K)}= M^2$. Differentiating the criterion and the constraint w.r.t. $\phi$, we get that such a critical point $\phi_o$ is given by the equation:
\[
\left\langle b_2b_2^*[\phi_o],\psi\right\rangle_{L^2(K)}-\left\langle b_2[{\be}],\psi\right\rangle_{L^2(K)} = -\lambda  \, \langle \phi_o,\psi\rangle_{L^2(K)} \, , \]
with $\gamma\in \mathbb{R}$ and for all $\psi \in L^2(K)$.
It follows that: 
\begin{equation}\label{cpe_ps}
\left\langle b_2 \, b_2^*[{\phi_o}]+ \lambda \, \phi_o,\psi\right\rangle_{L^2(K)}=\langle b_2({\be}),\psi\rangle_{L^2(K)},\,\,\,\hbox{for all } \psi \in L^2(K) \, ,
\end{equation}
\noindent whence: 
\begin{equation}\label{CPEmo}
  \tag{CPE$_{mo}$}
  b_2 \, b^*_{2}[\phi_o] +\lambda   \, \phi_o = b_2[{\be}] \, .
\end{equation}
\noindent We have $\lambda > 0$. Indeed, taking $\psi=\phi_o$ in Equation \eqref{cpe_ps}, we get 
$$-\lambda  \,  \langle\phi_o,\phi_o\rangle_{L^2(K)}=-\lambda  \,  \Vert \phi_o\Vert_{L^2(K)}^2=\langle b_2^*(\phi_o)-{\be},b_2^*(\phi_o)\rangle_{L^2(S,\mathbb{R}^2)}.$$
\noindent By \cite[Thm 5.2]{Brezis}, for any ${\bf g}\in L^2(S,\RR^2)$,
$$\Re\langle b_2^*(\phi_o)-{\be},b_2^*(\phi_o)-{\bf g}\rangle_{L^2(S,\mathbb{R}^2)}\leq 0.$$
\noindent Taking ${\bf g}=0$, we obtain that $\langle b_2^*(\phi_o)-{\be},b_2^*(\phi_o)\rangle_{L^2(S,\mathbb{R}^2)}\leq 0$. Thus, $-\lambda  \, \leq 0$. As $\lambda=0$ if and only if ${\be}\in \hbox{Ran}\,b_2^*$, we conclude that $-\lambda  < 0$.\\

\noindent Assume that ${\be}\not\in \hbox{Ran}\,b_2^*$. Then there exists a unique $\lambda > 0$ such that \eqref{CPEmo} holds true and $\Vert\phi_o\Vert_{L^2(K)}=M$. In this case,
$$\lambda\, \Vert \phi_o\Vert_{L^2(K)}^2= \lambda\, M^2=- \langle b_2^*(\phi_o)-{\be},b_2^*(\phi_o)\rangle_{L^2(S,\mathbb{R}^2)} \, .$$
\noindent From now on,  when there is no ambiguity, we will write $f^\prime = \partial_x \, f $ for the derivative of real-valued differentiable functions $f$ defined on $\RR$ or on intervals of $\RR$.\\

\noindent For $f\in L^2(K)$, at $y =h >0$, we have that:
\begin{equation*} 
b_2 \, b_2^*[f] = \left[
  - P_h \star  (\chi_S \, (P_h \star \widetilde{f})^\prime)
  + Q_h \star ( \chi_S \, (Q_h^\prime \star \widetilde{f}))
  \right]^\prime_{|_K}\end{equation*}
\[= \left[
  - P_h^\prime \star  (\chi_S \, (P^\prime_h \star \widetilde{f}))
  + Q_h^\prime \star(\chi_S \, (Q^\prime_h \star \widetilde{f}))
  \right]_{|_K}\, .\]
\noindent The critical point equation \eqref{CPEmo} can be written as:
\begin{equation*}
 \left[ - P_h^\prime \star  (\chi_S \, (P^\prime_h \star \widetilde{\phi_o}))
  + Q_h^\prime \star(\chi_S \, (Q^\prime_h \star \widetilde{\phi_o}))
  \right]_{|_K}+\lambda\phi_o=-P_h'\star(\widetilde{e_1}-\cH\widetilde{e_2})_{|_K},
  \end{equation*}
\noindent where $(e_1,e_2)={\be}$.\\

\noindent Observe that, 
if the  function $\eta \in L^2(K)$ represents the error on the measurements of $b_2[{\bm}]$, then 
\begin{eqnarray*}
\left\vert \langle b_2[{\bm}]+\eta \, , \, \phi_{i} \rangle_{L^2(K)} - \langle m_i \rangle\right\vert&=&\left\vert\langle{\bm} \, , \, b_2^{*}[\phi_{i}] - {\be}_i \rangle_{L^2(S,\RR^2)}+\langle\eta,\phi_{i}\rangle_{L^2(K)}\right\vert\\
&\leq&\Vert {\bm}\,\Vert_{L^2(S,\RR^2)}\Vert b_2^*[\phi_{i}]-{\be}_i\Vert_{L^2(S,\RR^2)} +\Vert \eta\Vert_{L^2(K)}\, \Vert \phi_{i}\Vert_{L^2(K)}\\
&\leq&\Vert {\bm}\,\Vert_{L^2(S,\RR^2)}\Vert b_2^*[\phi_{i}]-{\be}_i\Vert_{L^2(S,\RR^2)} +M\,\Vert \eta\Vert_{L^2(K)}. 
\end{eqnarray*}
\noindent For a choice of $M>0$, the solution of \eqref{BEP} guarantees that $\Vert b_2^*[\phi_{i}]-{\be}_i\Vert_{L^2(S,\RR^2)}$ is minimal and the term $\Vert \eta\Vert_{L^2(K)}$ is controlled by $M=\Vert\phi_{i}\Vert_{L^2(K)}$. The challenge of the numerical implementations is to choose a right $M>0$ such that $\phi_{i}$ satisfies the two conditions mentioned above. We will observe in numerical implementations that the solution $\phi_{i}$ can oscillate at the endpoints. 
However, since $\text{Ran} \, b_2 \subset W^{1,2}(K)$, then  \eqref{CPEmo} is to the effect that $\phi_{i} \in W^{1,2}(K) \subset C(\overline{K})$.

\subsection{BEP in $W_0^{1,2}(K)$}

\noindent An alternative is to formulate the BEP in a space of functions such that these oscillations are controlled.
It thus seems natural to search for a solution in the Sobolev space $W_0^{1,2}(K)$.

\noindent We denote by $b^*_{2|_W}$ the restriction of $b_2^*$ to $W_{0}^{1,2}(K)$. For $\phi\in W_{0}^{1,2}(K)$,
\[
  b_2^*[\phi]=b^*_{2|_W}[\phi]=((P_h\star \widetilde{\phi})'_{|_S},(Q_h\star\widetilde{\phi})'_{|_S}) 
  =((P'_h\star \widetilde{\phi})_{|_S},(Q'_h\star\widetilde{\phi})_{|_S})
  =((P_h\star\widetilde{\phi'})_{|_S},(Q_h\star\widetilde{\phi'})_{|_S}),
  \]
\noindent where $(\widetilde{\phi})^{\prime}=\widetilde{\phi^\prime}=\phi^\prime\vee 0$. Notice that $\widetilde{\phi}\in W^{1,2}(\RR)$ only when $\phi\in W^{1,2}_0(K)$ (see  \cite{Brezis}). 

\noindent Let  ${\be} \in L^2(S, \RR^2)$ and $M\geq 0$. The problem can thus be stated as the one of finding  $\phi_{o} \in W_{0}^{1,2}(K)$,  $\Vert \phi_{o} \Vert_{W_0^{1,2}(K)}\leq M$ such that:
\[
\Vert  b^*_{2|_W}[{\phi_{o}}] -{\be}\Vert_{L^2(S,\RR^2)}= \min_{\phi\in W_{0}^{1,2}(K)\atop{
\Vert \phi\Vert_{W_0^{1,2}(K)}\leq M}}\Vert b^*_{2|_W}[{\phi}] -{\be}\Vert_{L^2(S,\RR^2)} \, .\label{BEP1}\tag{BEP$_{mo,W}$}
\]

\subsubsection{Well-posedness}

\begin{proposition}\label{prop:wpBEP1}
  Let ${\be} \in L^2(S, \RR^2)$ and $M\geq 0$. Problem \eqref{BEP1} 
  admits a unique solution $\phi \in W^{1,2}_0(K)$.
  If
 $ {\be}  \not \in  b^*_{2|_W}[W_0^{1,2}(K)]$, then $\phi_o$ saturates the constraint:~$\Vert \phi_o \Vert_{W_0^{1,2}(K)}= M$.
\end{proposition}

\beginpf 
We apply again \cite[Lemma 2.1]{CP} in order to get  existence and uniqueness of a solution $\phi_o\in W^{1,2}_0(K)$  to \eqref{BEP1} saturating the constraint if $ {\be}  \not \in  b^*_{2|_W}[W_0^{1,2}(K)]$. Yet, because $b_2^*$ is injective (from Proposition \ref{prop:b2*inj}),
  uniqueness still holds when ${\be} \in b^*_{2|_W}[W_0^{1,2}(K)]$. Note that in this case, the constraint is saturated also  if ${\be} =b_{2}^*[\phi]$ for $\phi \in {L^2(K)}$ such that $\Vert \phi\Vert_{W_0^{1,2}(K)}\geq M$.
\endpf
\begin{remark}
  \noindent Let $\psi\in L^2(K)$ with zero mean on $K$. There exists a unique $\phi\in W_0^{1,2}(K)$ such that $\phi^\prime=\psi$ 
  (and $\Vert\phi\Vert_{W_0^{1,2}(K)} = \Vert\phi^\prime\Vert_{L^2(K)}$). 
  Then, \eqref{BEP1} is equivalent to the following problem: let  ${\be} \in L^2(S, \RR^2)$ and $M\geq 0$, find  $\psi_{o} \in {L^2(K)}$ such that $\int_K \psi_{o}(x) dx=0$, $\Vert \psi_{o} \Vert_{L^2(K)}\leq M$ that satisfies
\begin{equation}\label{a2*:min}
\Vert a_2^*[\widetilde{\psi_{o}}]-{\be}\Vert_{L^2(S,\RR^2)}= \min_{\psi\in L^2(K)\atop{\int_K \psi(x) dx=0\\
,\,\,\,\Vert \psi\Vert_{L^2(K)}\leq M}}\Vert a_2^*[\widetilde{\psi}] -{\be}\Vert_{L^2(S,\RR^2)} \, .
\end{equation}
\noindent Proposition \ref{prop:wpBEP1} is to the effect that a unique solution $\phi_o$ to \eqref{BEP1} does exist, whence also a unique solution $\psi_o$ to \eqref{a2*:min}, which is such that $\phi_o^\prime=\psi_o$.
\end{remark}
\noindent We also have the following density result.
\begin{lemma}
  \label{lem:denss}
    $b_2^*\left[W_{0}^{1,2}(K)\right]$ is dense in $L^2(S,\RR^2)$.
  \end{lemma}
  \beginpf
  \noindent First, observe that $W_{0}^{1,2}(K)$ is dense in $L^2(K)$ since the set $C_c^\infty(K)$ of $C^\infty(K)$ functions with compact support $\subset K$ is dense in $L^2(K)$ (see \cite[Cor. 4.23]{Brezis}). As $W^{1,2}(K)\subset L^2(K)$, we use the continuity of $b_2^*$ to get:
    \[
 b_2^*\left[{W_{0}^{1,2}(K)}\right] \subset  b_2^*(L^2(K))=b_2^*\left[\overline{W_{0}^{1,2}(K)}^{L^2(K)}\right] \subset
    \overline{b_2^*\left[W_{0}^{1,2}(K)\right]}^{L^2(S,\mathbb{R}^2)} \, ,
    \]
    Recall that ${b_2^*[L^2(K)]}$ is dense in $L^2(S, \RR^2)$ (for the $L^2(S, \RR^2)$ norm), 
    from Corollary \ref{cor:denserange}. We therefore conclude to the density of $b_2^*[W_{0}^{1,2}(K)]$ in $L^2(S, \RR^2)$.
\endpf\\      

\noindent Yet, Lemma \ref{lem:denss} 
implies that without the norm constraint, Problem \eqref{BEP1} is ill-posed. The minimum is in this case an infimum and is equal to 0 while minimizing sequences diverge.

\subsubsection{A constructive solution} 
\label{sec:csolw}
\noindent Observe that ${b_2}=({b_2^*}_{|_W})^*$, since ${b_2^*}_{|_W}=a_2^*(\phi')$, where
$\phi'  $ denotes the derivation from $W_0^{1,2}(K)$ onto the space of $L^2(K)$-functions with zero mean.
Thus, for $f\in W_{0}^{1,2}(K)$, we get
\begin{equation*}
b_2 \, b^*_{2|_W}[f] 
= \left[
  - P_h^\prime \star (\chi_S \, P_h \star \widetilde{f^\prime})
  + Q_h^\prime  \star (\chi_S \, Q_h \star \widetilde{f^\prime})
  \right]_{|_K}.
\end{equation*}
\noindent Assume that $ {\be}  \not \in  b^*_{2|_W}[W_0^{1,2}(K)]$.
From Proposition \ref{prop:wpBEP1}, we obtain a  critical point equation  somewhat different from \eqref{CPEmo} for the solution $\phi_o$ in this case (see also Section \ref{sec:csol}): 
\[
\left\langle b_2 b^*_{2|_W}[{\phi_o}],\psi\right\rangle_{L^2(K)}-\left\langle b_2[{\be}],\psi\right\rangle_{L^2(K)} = \gamma \, \langle \phi_o',\psi'\rangle_{L^2(K)} = -\gamma \, \langle \phi_o'',\psi\rangle_{L^2(K)} \, , \]
\noindent with $\gamma<0$ and $\forall \psi \in W_0^{1,2}(K)$. The solution $\phi_o\in W_0^{1,2}(K)$ of \eqref{BEP1} is thus given by  
\begin{equation}\label{CPEmoW}
  \tag{CPE$_{mo, W}$} %
b_2 \, b^*_{2|_W}[\phi_o] + \gamma \phi_o'' = b_2[{\be}] \, ,
\end{equation}
with $\gamma <0$, such that $\Vert \phi_o \Vert_{W^{1,2}(K)}=\Vert  \phi_o' \Vert_{L^2(K)}= M$ if ${\be} \not \in \hbox{Ran}\,b^*_{2|_W}$. The equation \ref{CPEmoW} is considered in the weak sense in $L^2(K)$. 
Thus, 
$$\left[
  - P_h^\prime \star  (\chi_S \, (P_h \star \widetilde{\phi'_o}))
  + Q_h^\prime \star(\chi_S \, (Q_h \star \widetilde{\phi'_o}))
  \right]_{|_K}+\gamma\phi_o''=-P_h'\star(\widetilde{e_1}-\cH\widetilde{e_2})_{|_K}. $$

\subsubsection{An equivalent problem in $H^2$}

\noindent Let $\mathbb{C}_+=\{x+iy\in\mathbb{C}: x\in\mathbb{R},y>0\}$ and $H^2(\mathbb{C}_+)$ denotes the Hardy space of the upper half-plane defined in Section \ref{sec:hardy}. As a function $F\in H^2(\CC_+)$ can be identified to its boundary function $F^*\in H^2(\RR)$ so that $F_{|_{S\times\{0\}}}$ can be defined as $F^*_{|_{S}}$.\\

\noindent Given a function $g\in L^2(S,\mathbb{C})$, for $\varepsilon>0$, there is a function $F\in H^2(\mathbb{C}_+)$ such that $\Vert g-F^*_{|_{S}}\Vert_{L^2(S,\mathbb{C})}<\varepsilon$. Indeed, if $f\in \left(H^2(\RR)_{|_S}\right)^{\perp}$, then for $g\in H^2(\RR)$, 
$$\langle f,g_{|_S}\rangle_{L^2(S,\CC)}=\langle \widetilde{f},g\rangle_{L^2(\RR,\CC)}=0. $$
\noindent Taking $g(t)=\frac{1}{t-z}$ with $z\in \CC$ such that $\Im z>0$, it follows that 
$$\int_{-\infty}^{+\infty}\frac{\widetilde{f}(t)}{t-\overline{z}}=0,\quad z\in\CC_+, $$
\noindent and $\widetilde{f}\in H^2_-(\RR)$. By the boundary uniqueness Theorem, $f=0$ on $S$ so $H^2(\RR)_{|_S}$ is dense in $L^2(S,\mathbb{C})$.\\ 

\noindent For $M \geq 0$, let
\[
\mathcal{B}_M=\{F\in H^2(\mathbb{C}_+):\text{supp} \, \Re F^{*} \subset K\,\,\text{ and }\,\,\,\Vert\Re F^{*} \Vert_{L^2(K)}\leq M \}\subset H^2(\mathbb{C}_+) \, .
\]
\begin{proposition}\label{bep:eqbis}
  Let 
  $g\in L^2(S,\CC)$. 
The problem \eqref{BEP1} is equivalent to the following bounded extremal problem (with $L^2(S,\CC)$ norm):

\[
\min_{F\in \mathcal{B}_M}\Vert F(.+ih)_{|_S}-g\Vert_{L^2(S,\CC)},\label{BEP2bis}\tag{BEP$_{\CC_+}$bis} \, .
\]
\end{proposition}
\beginpf Let ${\be}\in L^2(S,\RR^2)$ be such that $J(g)={\be}$. By Proposition \ref{prop:wpBEP1}, there exists a function ${\phi_o}\in W_0^{1,2}(K)$ such that 
$$\Vert {b_2^*}_{|_W}[\phi_o]-{\be}\Vert_{L^2(S,\RR^2)}=\min_{\phi\in W_{0}^{1,2}(K)\atop{
\Vert \phi\Vert_{W_0^{1,2}(K)}\leq M}}\Vert b^*_{2|_W}[{\phi}] -{\be}\Vert_{L^2(S,\RR^2)}\, . $$
\noindent Using the isomorphism $J$, we get for any $\phi\in W_0^{1,2}(K)$,
\[\Vert {b_2^*}_{|_W}[\phi] -{\be}\Vert_{(L^2(S, \RR^2)}=\Vert J({B_2^{*}}_{|_W}[\phi])-{\be}\Vert_{L^2(S,\RR^2)}
= \Vert {B_2^{*}}_{|_W}[\phi] -J^{-1}{\be}\Vert_{L^2(S,\CC)} \,  ,
\]
\noindent where ${B_2^{*}}_{|_W}[\phi]=(P_h\star(I+i\mathcal{H})\widetilde{\psi})_{|_S}$ and $\psi$ is the (unique) function of $L^2(K)$ such that $\phi'=\psi$. Let $F_o \in {\mathcal{B}}_M$ defined by $P_y\star (I+i\mathcal{H})\widetilde{\psi_o}=F_o$ with $\phi_o'=\psi_o$. It follows that $F_o(.+ih)={B_2^*}_{|_W}[\phi_o]=(P_h\star(I+i\mathcal{H})\widetilde{\psi_o})$ and ${(\Re \, F_o^*)}_{|_K}=\psi_o$. Then, 
\begin{eqnarray*}
\Vert F_o(.+ih)_{|_S}-g\Vert_{L^2(S,\CC)}&=&\Vert {B_2^{*}}_{|_W}[\phi_o] -J^{-1}{\be}\Vert_{L^2(S,\CC)}\\
&=&\min_{\phi\in W_{0}^{1,2}(K)\atop{
\Vert \phi\Vert_{W^{1,2}(K)}\leq M}}\Vert B^*_{2|_W}[{\phi}] -J^{-1}{\be}\Vert_{L^2(S,\RR^2)}\\
&=&\min_{\phi\in W_{0}^{1,2}(K)\atop{
\Vert \phi\Vert_{W^{1,2}(K)}\leq M}}\Vert (P_h\star(I+i\mathcal{H})\widetilde{\phi'})_{|_S} -J^{-1}{\be}\Vert_{L^2(S,\RR^2)}\\
&=&\min_{\psi\in L^2(K)\atop{\Vert \psi\Vert_{L^2(K)}\leq M}}\Vert (P_h\star(I+i\mathcal{H})\widetilde{\psi})_{|_S} -J^{-1}{\be}\Vert_{L^2(S,\RR^2)}\\
&=&\min_{F\in \mathcal{B}_M}\Vert F(.+ih)_{|_S}-g\Vert_{L^2(S,\CC)},
\end{eqnarray*}
\noindent where we use that $P_y\star (I+i\cH)\widetilde{\psi}\in {\mathcal{B}}_{M}$ for any $\psi\in L^2(K)$. 
\endpf

\section{Algorithms and numerical computations}
\label{sec:comput}

\subsection{Construction of the algorithms}
\subsubsection{Bases functions}

\noindent   Let $S = (-s,s)$, $K=(-q,q)$, $s$, $q \in \RR_+$. For $n\in \ZZ$, let $g_n(x)=\exp \, (i \, n \, \pi \, x / q)$, for $x \in K$.
They are eigenfunctions of the Laplacian on $K$:
\[
g_n''= - \left[\frac{n \pi}{q} \right]^2 \, g_n= - \mu_n\,  g_n\, , \quad\hbox{where}\quad \mu_n= \left[\frac{n\pi}{q} \right]^2 >0 \, .
\]
\noindent Up to multiplication by the constant factor $\frac{1}{({2 q})^{1/2}}$ for normalization, the family of functions $\left( g_n\right)_n$,  $n \in \ZZ$, is the Hilbertian Fourier basis of $L^2(K, \CC)$ 
(see \cite[Sec. 8.6]{Brezis}).
\noindent We  use it to expand and compute the solutions to the above bounded extremal problems,  \eqref{BEP} in $L^2(K)$ and \eqref{BEP1} in $W_{0}^{1,2}(K)$. 
\noindent {Indeed, functions 
$ \phi $ 
in $L^2(K)$ can be expanded in Fourier series on the basis $\left({g_n}\right)_n$,  $n \in \ZZ$: 
\begin{equation*}
 \phi 
          ~=~ \sum_{n \in \mathbb{Z}} c_n \, g_n \, , 
\end{equation*}
with $c_n  \in l^2(\mathbb{Z})$, $c_{-n} = \bar{c}_n$, $c_0 \in \RR$ (since $\phi$ is real-valued). 
\noindent Moreover, $ \phi $ belongs to $W^{1,2}(K)$  if, and only if, the coefficients of its expansion on $\left({g_n}\right)_n$,  $n \in \ZZ$, are such that  
$n \, c_n  \in l^2(\mathbb{Z})$} and $c_n = - c_{-n} \in i \, \RR$.

  \noindent Another family of appropriate functions in $L^2(K)$ is made of piecewise constant functions on small intervals (say $K_l$, $l = 1, \cdots, L$) covering $K=(-q,q)$ (so called $P_0$ finite elements), see also \cite{BEP-moments}. We use it for the computations of the solutions to the associated forward problem (Section \ref{sec:fp}).
             
          \subsubsection{Operators, matrices in $L^2(K)$ and $W_0^{1,2}(K)$}\label{ssec:cpel2}
         
\noindent             For all ${\phi} \in L^ 2(K)$, we use the expression
  \[
    {b_2^*}  [{\phi}]= \left((\partial_x \, P_y) \star \widetilde{\phi}, (\partial_x \, Q_y) \star\widetilde{\phi}\right)_{|_{S\times\{h\}}}
    = \left( P_y' \star \widetilde{\phi}, Q_y' \star\widetilde{\phi}\right)_{|_{S\times\{h\}}} \, .    \]

    \noindent Solutions to \eqref{BEP} in $L^2(K)$ (actually in $W^{1,2}(K) \subset L^2(K)$) are provided by the critical point equation \eqref{CPEmo}.

     \noindent    For moments estimation, we consider the specific functions ${\be} ={\be}_1 \, , {\be}_2$, with ${\be}_1 = (\chi_S,0)$ for $\langle m_1 \rangle$,  ${\be}_2 = (0,\chi_S)$ for $\langle m_2 \rangle$. 
    From Lemma \ref{lem:ph} and the relation \eqref{eg1}, it holds that
    ${\be}_i \not \in \hbox{Ran}\,b_2^*$. Hence, the constraint in the BEP is saturated by the approximant, following Propositions \ref{prop:wpBEP} and \ref{prop:wpBEP1}.

\noindent Hence, using the Fourier basis  $(g_n)$ in order to express solutions to (CPE$_{mo}$) as $\phi_o = \sum_{n \in \mathbb{Z}} c_n \, g_n$, we get that, $\forall k \in \ZZ$:
\begin{equation}
  \label{eq:cn}
\sum_{n \in \mathbb{Z}} c_n \, \left\langle {b_2^*}[{g_n}],{b_2^*}[g_k]\right\rangle_{L^2(S,\mathbb{R}^2)}+\lambda \, c_k
=\langle {\be},{b_2^*}[g_k]\rangle_{L^2(S,\RR^2)} \, ,
\end{equation}
for $\lambda > 0$ such that $\Vert\phi_o\Vert_{L^2(K)}^2=\sum_{n \in \mathbb{Z}} \left|c_n\right|^2 =M^2$.

\noindent For $t_1 \, , \ t_2 \in K$, $x \in S$, let:
  \[
  I(t_1,t_2,x) = P_h'(x-t_1) \,P_h'(x-t_2)+Q_h'(x-t_1) \, Q_h'(x-t_2) \, ,\]
and
\begin{equation*}
k(t_1,t_2)~=~\int_{-s}^{s}
I(t_1,t_2,x) \, dx \, .
\end{equation*}
\noindent Following the expressions of the operator $b_2^*$ in Section \ref{sec:mo}, we obtain:
\begin{lemma}
  \label{lem:Fourierb2*}
\begin{equation}\label{Fourierb2*}
  \left\langle {b_2^*}[{g_n}],{b_2^*}[g_k]\right\rangle_{L^2(S,\RR^2)}  ~=~  \frac{1}{\pi^2}
 \int_{-q}^{q}\int_{-q}^{q}
   k(t_1,t_2) \, e^{i \frac{\pi}{q}(n t_1 + k t_2)} \, dt_1 \, dt_2 \, .  
 \end{equation}\end{lemma}
These are the Fourier coefficients of  $\chi_K(t_1) \, \chi_K(t_2) \, k(t_1,t_2)$.\\

\beginpf \noindent Observe that 
\begin{eqnarray*}
\left\langle {b_2^*}[\phi],{b_2^*}[\psi]\right\rangle_{L^2(S,\mathbb{R}^2)}&=&\left\langle P_h' \star \widetilde{\phi},P_h' \star \widetilde{\psi}\right\rangle_{L^2(S)}+\left\langle Q_h' \star \widetilde{\phi},Q_h' \star \widetilde{\psi}\right\rangle_{L^2(S)}\\
&=&\int_{-s}^s P_h' \star \widetilde{\phi}(x)P_h' \star \widetilde{\psi}(x) dx+\int_{-s}^s Q_h' \star \widetilde{\phi}(x)Q_h' \star \widetilde{\psi}(x) dx
\end{eqnarray*}

\[
=\int_{-s}^s \int_{-q}^q\int_{-q}^q\left(P_h'(x-t_1)P_h'(x-t_2)+Q'_h(x-t_1)Q'_h(x-t_2)\right)\phi(t_1)\psi(t_2) dt_1dt_2dx \, .
\]
We have, for $\phi,~\psi \in L^2 (K)$,
 \[
\langle b_2^*[\phi],b_2^*[\psi]\rangle_{L^2(S, \RR^2)}
= \frac{1}{\pi^2}
 \int_{-s}^{s}\int_{-q}^{q}\int_{-q}^{q}
 I(t_1,t_2,x) \, \phi(t_1) \, \psi(t_2) \, dt_1 \,dt_2 \,dx  \, .
\]
We can interchange the integrals, which leads to
\begin{equation*}
\langle b_2^*[\phi],b_2^*[\psi]\rangle_{L^2(S, \RR^2)} ~=~  \frac{1}{\pi^2}
 \int_{-q}^{q}\int_{-q}^{q}
   k(t_1,t_2) \, \phi(t_1) \, \psi(t_2) \,  dt_1 \, dt_2 \, . 
\end{equation*}
Taking $\phi =g_n$ and $\psi = g_k$ leads to the expression \eqref{Fourierb2*}. \endpf\\
\noindent This allows to   compute the left-hand side of \eqref{CPEmo} and of \eqref{eq:cn}.\\

\noindent  Similarly, solutions to \eqref{BEP1} in $W_{0}^{1,2}(K)$, are furnished by \eqref{CPEmoW}.
Notice that $    {b_2^*}_{|_W} [{\phi}]= {b_2^*}  [{\phi}]$ if ${\phi} \in W_{0}^{1,2}(K)$.
Using   $(g_n)$  in order to express the solutions to (CPE$_{mo,W}$) as $\phi_o = \sum_{n \in \mathbb{Z}} d_n \, g_n$,  we get that, $\forall k \in \ZZ$:
\begin{equation}
  \label{eq:dn}
\sum_{n \in \mathbb{Z}} d_n \, \left\langle {b_2^*}[{g_n}],{b_2^*}[g_k]\right\rangle_{L^2(S,\RR^2)}- \gamma \, \mu_k \, d_k=\langle {\be},{b_2^*}[g_k]\rangle_{L^2(S,\RR^2)} \, . 
\end{equation}
for $\gamma <0$ such that $\Vert\phi_o'\Vert_{L^2(K)}^2=\sum_{n \in \mathbb{Z}} \mu_n \, \left|d_n\right|^2 =M^2$, which has to be compared  with \eqref{eq:cn}.

\noindent Note that the above functions $ I(t_1,t_2,x)$ and $k(t_1,t_2) $ admit explicit expressions which can be used to more precisely compute the quantities 
$\langle b_2^*[\phi],b_2^*[\psi]\rangle_{L^2(S, \RR^2)}$ in Lemma \ref{lem:Fourierb2*}.

\subsubsection{Forward problems, $b_2[{\be}_i]$, $b_2[{\bm}]$}
\label{sec:fp}

\noindent Concerning the computation of the right-hand sides $\left\langle {\be},{b_2^*}[g_k]\right\rangle_{L^2(S,\RR^2)}$ of the critical point equations, following the definition of the operator $b_2$ in Section \ref{sec:mo},
  \[b_2[{\bm}]= - \left(P'_h \star \widetilde{m_1} -Q'_h\star \widetilde{m_2}\right)_{|_{K}} \, , \ \ \forall {\bm}=(m_1,m_2) \in L^2(S, \mathbb{R}^2)\, ,  \]
we first compute, from the above definition of ${\be}_i$ for $i = 1,2$:
\[
\left\langle {\be}_1,{b_2^*}[g_k]\right\rangle_{L^2(S, \RR^2)} = \left\langle b_2[{\be}_1],g_k\right\rangle_{L^2(K)}
\]
\[= \int_{-s}^{s}\int_{-q}^{q} P_h'(x-t) \, e^{i k \pi \frac{t}{q}} \, dt \, dx = \int_{-q}^{q} \left(P_h(s-t) -P_h(s+t)\right) \, e^{i k \pi \frac{t}{q}} \, dt \, ,
\]
and
\[
\left\langle {\be}_2,{b_2^*}[g_k]\right\rangle_{L^2(S, \RR^2)} = \left\langle b_2[{\be}_2],g_k\right\rangle_{L^2(K)}
\]
\[= \int_{-s}^{s}\int_{-q}^{q} Q_h'(x-t) \, e^{i k \pi \frac{t}{q}} \, dt \, dx = \int_{-q}^{q} \left(Q_h(s-t) -Q_h(s+t)\right) \, e^{i k \pi \frac{t}{q}} \, dt \, .
  \]
\noindent   These expressions are computed by discretization on a family of piecewise constant functions on $K = \cup_l K_l$, summing up the contributions on $K_l$, although   they are    linked to the Fourier coefficients of $\chi_K \, P_h$ and $\chi_K \, Q_h$. 

\noindent   This is all we need in order to solve the linear systems \eqref{eq:cn} in $c_n$ and \eqref{eq:dn} in $d_n$, and to compute the solutions $\phi_{i}$ related to ${\be}_i$ and $M$.\\

\noindent In order to generate synthetic data, the components $m_i$ of magnetizations ${\bm} \in L^2(S, \RR^2)$ are also modelled as  finite linear combinations of  a family of piecewise constant functions on $S= \cup_l S_l$, $l= 1, \cdots L$. %

\noindent The quantity  $b_2[{\bm}]$ is computed from its Fourier series, while its Fourier coefficients (equal to $\left\langle b_2[{\bm}],g_k\right\rangle_{L^2(S)}$) are computed as those of $b_2[{\be}_i]$ above, summing up the contributions of the intervals $S_l$, $l= 1, \cdots L$.

\noindent We  finally compute (with a quadrature method):
\[\langle b_2[{\bm}] \, , \, \phi_{i} \rangle_{L^2(K)}= \sum_{n \in \mathbb{Z}} c_n \, \langle b_2[{\bm}] \, , \, g_n \rangle_{L^2(K)} \, ,
\]
which is expected to furnish the approximation $\langle m_i^{e} \rangle=\langle b_2[{\bm}] \, , \, \phi_{i} \rangle_{L^2(K)}$ to $\langle m_i \rangle$.
For error estimation,  we will compare $\langle m_i^{e} \rangle$ to the actual moments $\langle m_i \rangle$ ($i = 1,2$).

\subsection{Numerical illustrations}
\label{sec:numres}

\noindent In this section, we discuss preliminary numerical results concerning the computation of  $\langle m_i^{e} \rangle=\langle b_2[{\bm}] \, , \, \phi_{i} \rangle_{L^2(K)}$ in order to estimate the mean value $\langle m_i\rangle$ of each component of ${\bf m}$.

\noindent  Recall that the solutions to (CPE$_{mo}$) in $L^{2}(K)$, $\phi_i=\phi_o[{\be}_i]$, $i=1,2$, 
\noindent are given by:

\[ \left[ - P_h^\prime \star  (\chi_S \, (P^\prime_h \star \widetilde{\phi_i}))
  + Q_h^\prime \star(\chi_S \, (Q^\prime_h \star \widetilde{\phi_i}))
  \right]_{|_K}+\lambda\phi_i=\left\{
\begin{array}{l}
  -\left( P_h'\star \chi_S \right)_{|_K} \, i = 1 \, ,\\
  \left( P_h'\star \cH \, \chi_S \right)_{|_K}\, i = 2 \, ,
\end{array}
\right. \]
\noindent    for $\lambda >0$; we have $\Vert \phi_i \Vert_{L^{2}(K)}= M (\lambda)$.

\noindent  The solutions $\phi_i=\phi_o[{\be}_i]$, $i=1,2$, to (CPE$_{mo,W}$) in $W_0^{1,2}(K)$ are given by:
 
\[ \left[
  - P_h^\prime \star  (\chi_S \, (P_h \star \widetilde{\phi'_i}))
  + Q_h^\prime \star(\chi_S \, (Q_h \star \widetilde{\phi'_i})) \right]_{|_K}-\lambda\phi_i''=
\left\{
\begin{array}{l}
  -\left( P_h'\star \chi_S \right)_{|_K} \, i = 1 \, ,\\
  \left( P_h'\star \cH \, \chi_S \right)_{|_K}\, i = 2 \, ,
\end{array}
\right.
\]
with $\lambda = - \gamma >0$, and yet $\Vert \phi_i \Vert_{W_0^{1,2}(K)}= M (\lambda)$, in the weak sense in $L^2(K)$. 

\noindent The quantities (coefficients) $\left\langle {b_2^*}[{g_n}],{b_2^*}[g_k]\right\rangle_{L^2(S,\RR^2)}$ in \eqref{eq:cn} and \eqref{eq:dn}
are computed with {\it fft} (fast Fourier transform), following Lemma \ref{lem:Fourierb2*}, while for numerical purposes, the series will be truncated at some order $N\in \NN_*$.

\noindent Various situations can be examined, with respect to the distance $h$ between the magnetization support $\bar{S}$ contained in 
$\RR \times \{0 \}$ and the measurement set $K \subset \RR \times \{h \}$,  and to their length $|S|$, $|K|$.
In order to match with the physics of the model, we assume that
$S= (-s,s) $ and $K= (-q,q) $ are centered at the origin, and that
$|K| > |S| > h >0$, with (approximately) $|K| \simeq 1.5 \, |S|$, $|S| \simeq 20 \, h$. We will take $h=0.1$, $s=1$ and $q=1.5$.
The order of truncation will be $N=250$.

\noindent We will compare the estimates $\langle m_i^e\rangle$ of $\langle m_i\rangle$ for functions $\phi_i$ solutions to the BEP in $L^2$ or in $W_0^{1,2}$, associated to particular values of $\lambda$, for some magnetizations $\bm$.
The computations are performed with Matlab, R2017a.

\noindent Tables below 
 furnish the estimated moments $\langle m_i^{e} \rangle$ for $i\in\{1,2\}$, the corresponding values of the relative errors $\varepsilon_i$ and norm constraints $M$, $\Vert \phi_i\Vert_{L^2(K)}$ or $\Vert \phi_i\Vert_{W^{1,2}_0(K)}$. Figures are provided in the Appendix.

\noindent 
Below, for the computation of the solutions $\phi_i$  to the corresponding BEP, we choose $\lambda=\lambda_i$ such that $1 \leq -\log_{10}(\lambda_i) \leq 9$ is an integer and such that the relative error for moments estimation given by
 \[\varepsilon_i =\frac{\left\vert\langle m_i\rangle - \langle m_i^e  \rangle\right\vert}{\vert\langle m_i\rangle\vert} \, 
  ,\]
 \noindent is small enough 
 among this range for the corresponding data $\bm$, while still providing an acceptable value for the quantity $M(\lambda_i)$, say between 10 and 20, see Table \ref{table0}.
 Observe that this trade off is made possible by the fact that the data are available for the present simulations. Elsewhere, in general, $\lambda=\lambda_i$ has to be chosen in terms of the behaviour of the error $\Vert b^*_{2}[{\phi_i}] -{\be}_i\Vert_{L^2(S,\RR^2)}$ (and of the associated constraint $M(\lambda_i)$), independently of the unknown magnetization $\bm$.

\noindent  In all the numerical experiments, the net moment will be kept fixed to: $\langle m_1\rangle=-0.1$ and $\langle m_2\rangle=0.1$, whence: 
\[\varepsilon_i
= 10 \, \left\vert \langle m_i  \rangle- \langle m_i^e  \rangle\right\vert
= \left\vert(-1)^{i-1} + 10 \, \langle m_i^e  \rangle\right\vert \,
.\]

\subsubsection{Solutions $\phi_i$ to BEP}
\label{subs:phii}
\begin{center}
\begin{figure}[h]
  \begin{tabular}{ll}    \includegraphics[scale=0.45]{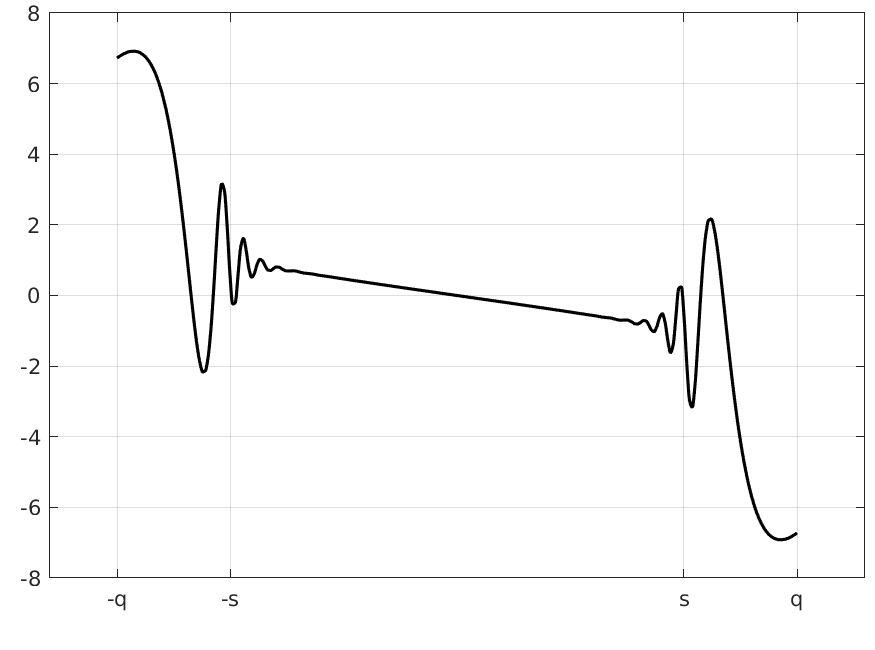} & \includegraphics[scale=0.45]{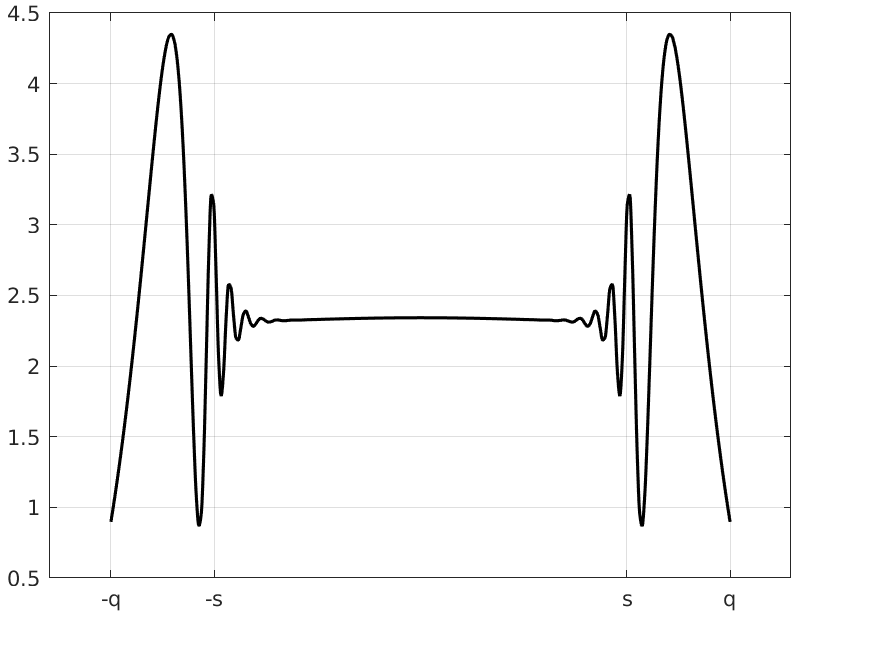} 
  \end{tabular}
  \caption{Solutions $\phi_i \in L^{2}(K)$ for $i=1,2$, $\lambda_i=10^{-3}$.}
    \label{figphi12-8-0}
\end{figure}
\end{center}

\begin{center}
\begin{figure}[h]
  \begin{tabular}{ll}    \includegraphics[scale=0.45]{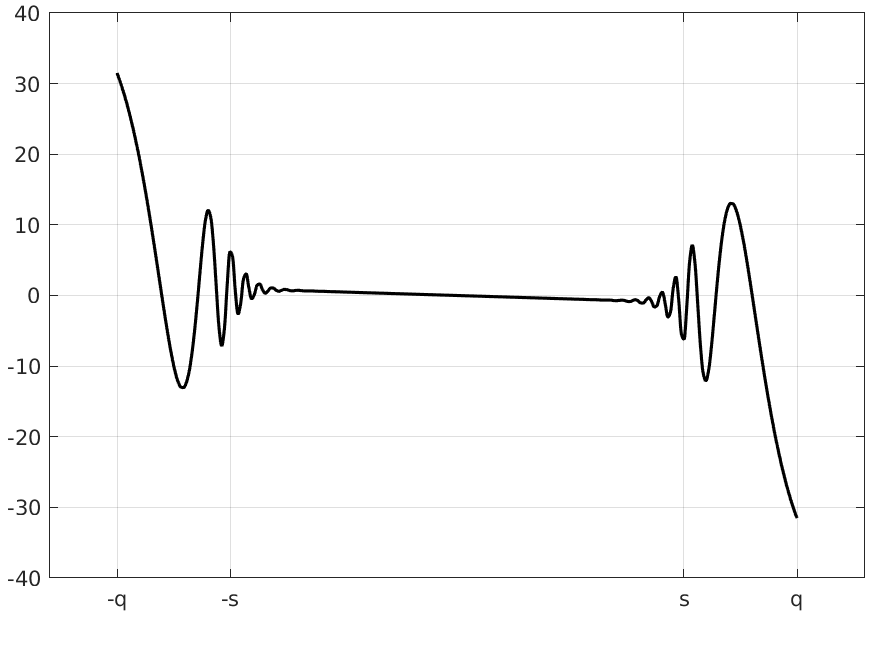}
    & \includegraphics[scale=0.45]{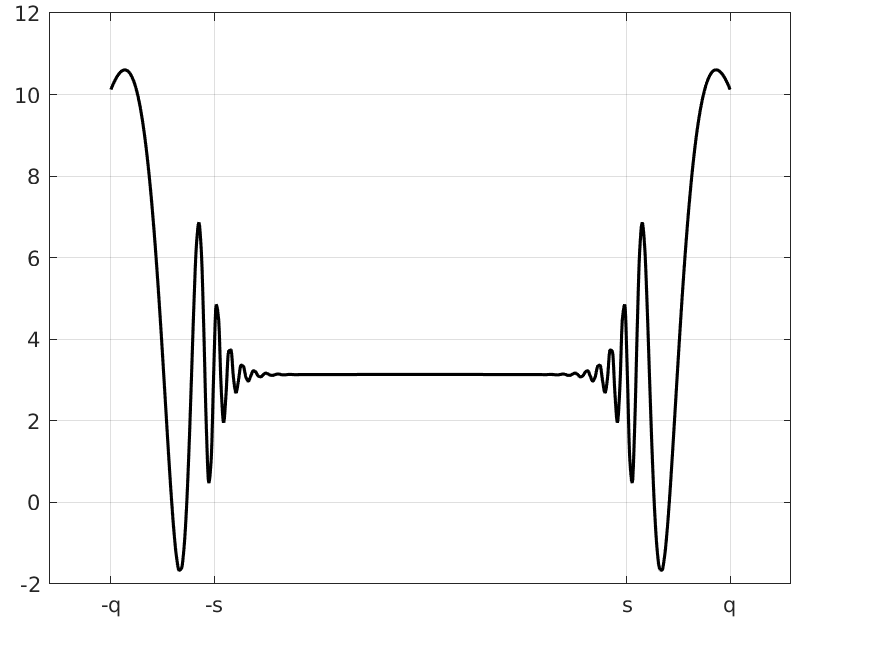}
      \end{tabular}
  \caption{Solutions $\phi_i \in L^{2}(K)$ for $i=1,2$, $\lambda_i=10^{-5}$.}
    \label{figphi12-8}
\end{figure}
\end{center}

\begin{center}
    \begin{figure}[h]
      \begin{tabular}{ll}   \includegraphics[scale=0.45]{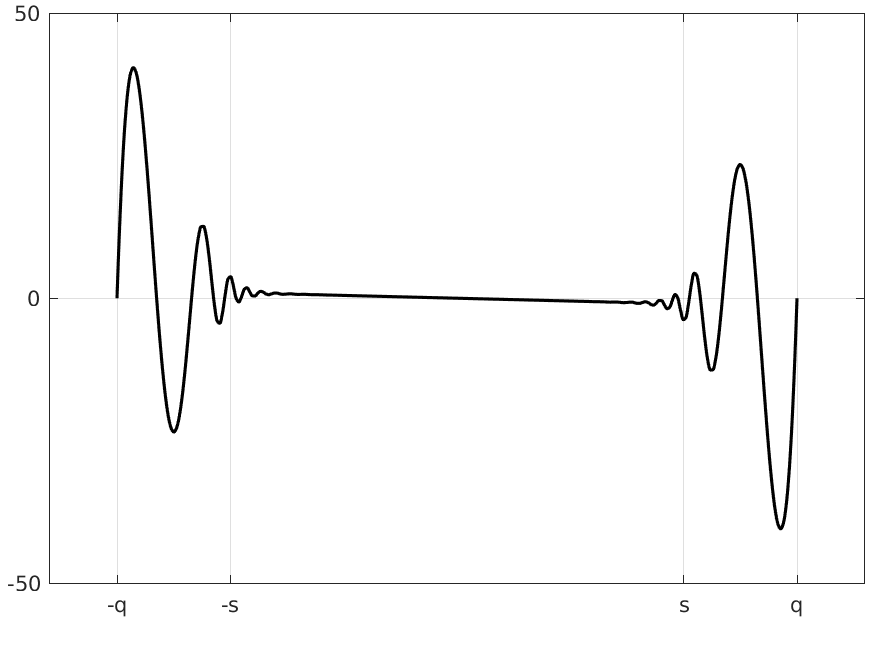}
        &  \includegraphics[scale=0.45]{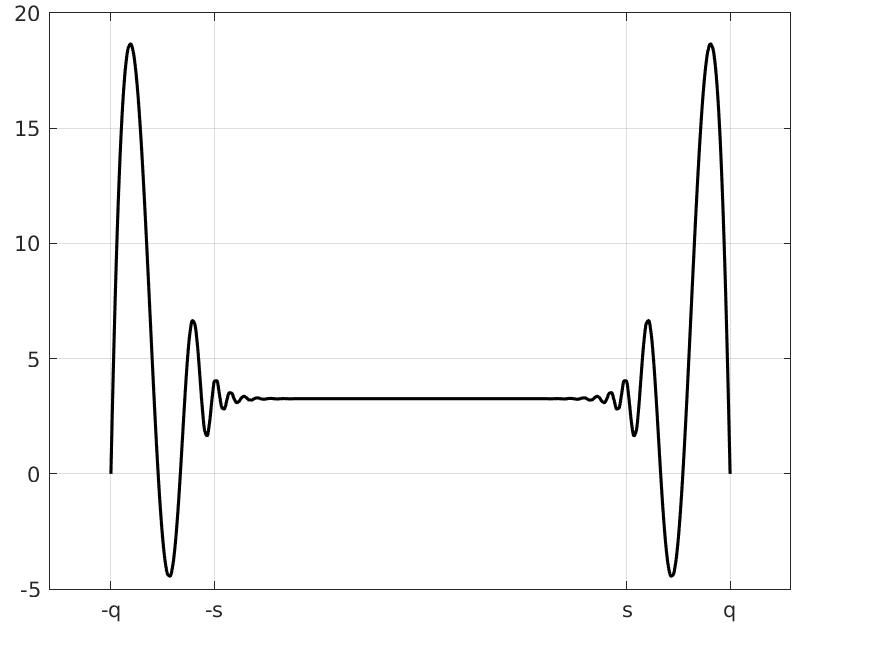}
            \end{tabular}
     \caption{Solutions $\phi_i \in W^{1,2}_0(K)$, $i=1,2$, $\lambda_i=10^{-8}$.} 
     \label{figphi12-bis}
    \end{figure}
\end{center}

\begin{center}
\begin{figure}[h]
  \begin{tabular}{ll}   \includegraphics[scale=0.45]{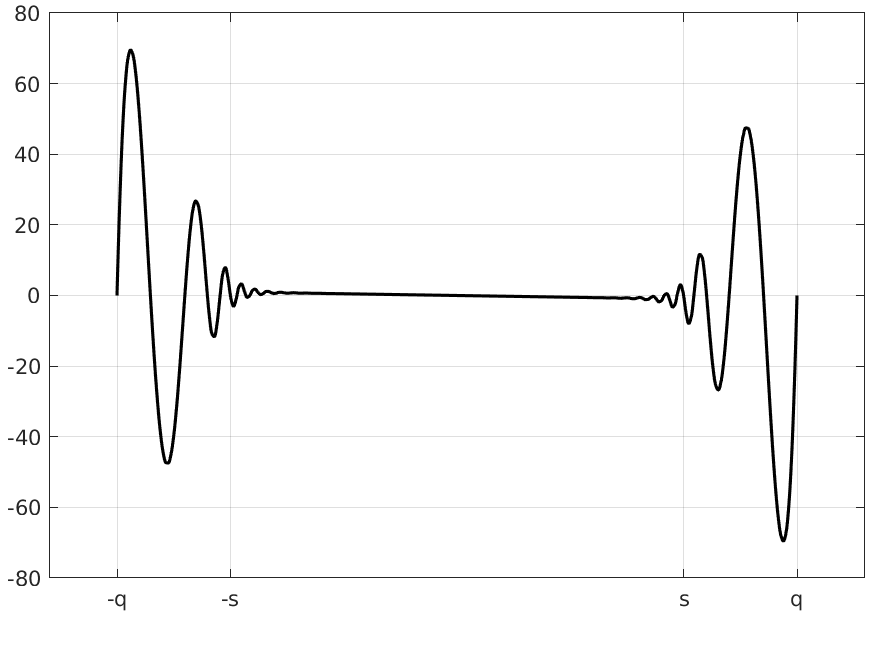}
    & 
    \includegraphics[scale=0.45]{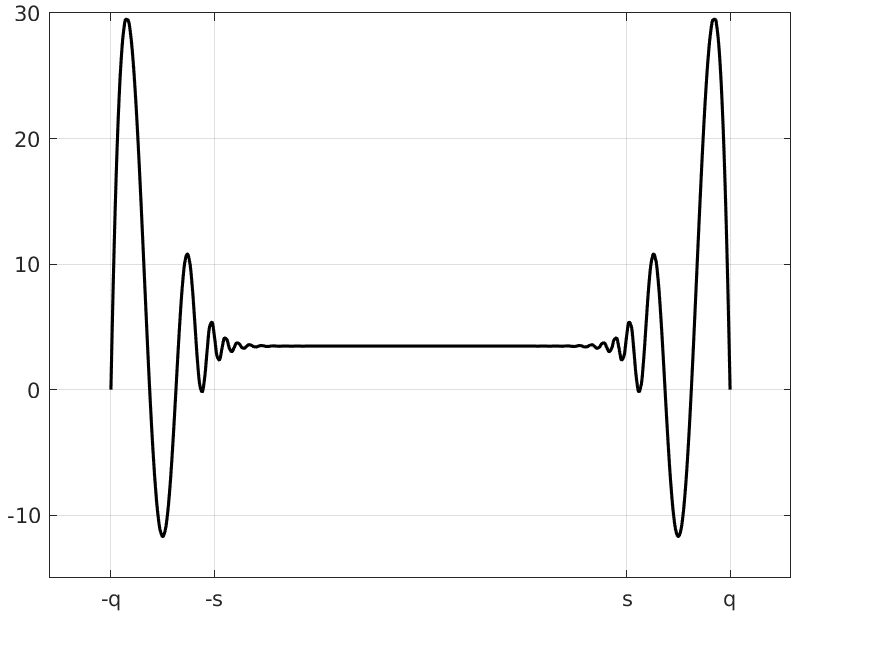}
     \end{tabular}
     \caption{Solutions $\phi_i \in W^{1,2}_0(K)$, $i=1,2$, $\lambda_i=10^{-9}$.} 
     \label{figphi12-bis-0}
\end{figure}
\end{center}

  \begin{table}[h]
    \centering
\begin{tabular}{|c|c|c|c|c|c|c|}
\hline
$\phi_i$ & $\lambda_i$ & $\Vert \phi_1\Vert$ & $\Vert \phi_2\Vert$ \\[0.07cm]
\hline
in $L^2(K)$ & $10^{-3}$& 4.8 & 4.4\\ \hline

 & $10^{-5}$& 14.4 & 8.2\\ \hline
in $W^{1,2}_0(K)$& $10^{-8}$ & 19.9 & 10.4\\ \hline
& $10^{-9}$ & 645.5 & 221.7\\ \hline
\end{tabular}
\caption{Values of $\lambda_i$, $M(\lambda_i)$= $\|\phi_i\|_{L^2(K)}$ or $\|\phi_i\|_{W^{1,2}_0(K)}$, $i=1,2$.}
\label{table0}
\end{table}
\noindent  The solutions $\phi_i\in L^{2}(K)$ for $i=1,2$ are plotted in Figure \ref{figphi12-8-0}, \ref{figphi12-8} ($\lambda_i=10^{-3}$ and $10^{-5}$). The solutions $\phi_i\in W^{1,2}_0(K)$ can be seen in Figures \ref{figphi12-bis},  \ref{figphi12-bis-0} ($\lambda_i=10^{-8}$ and $10^{-9}$). Their norms are given in Table \ref{table0}. 
 
\noindent  The solutions in $L^2(K)$ show quite many oscillations close to the boundary of $S$, less in $W^{1,2}_0(K)$,  as expected. Interestingly, these functions behave as affine or constant in a large interval contained in the interior of $S$. We also see that as $\lambda_i$ decreases, the corresponding constraint $M(\lambda_i)$ grows, which also increases the oscillating phenomenon, mainly in $L^2(K)$.

\subsubsection{Constant magnetizations} 
\label{subs:mwvls}

\begin{figure}[h]
          \begin{tabular}{lll}    \includegraphics[width=.3\textwidth]{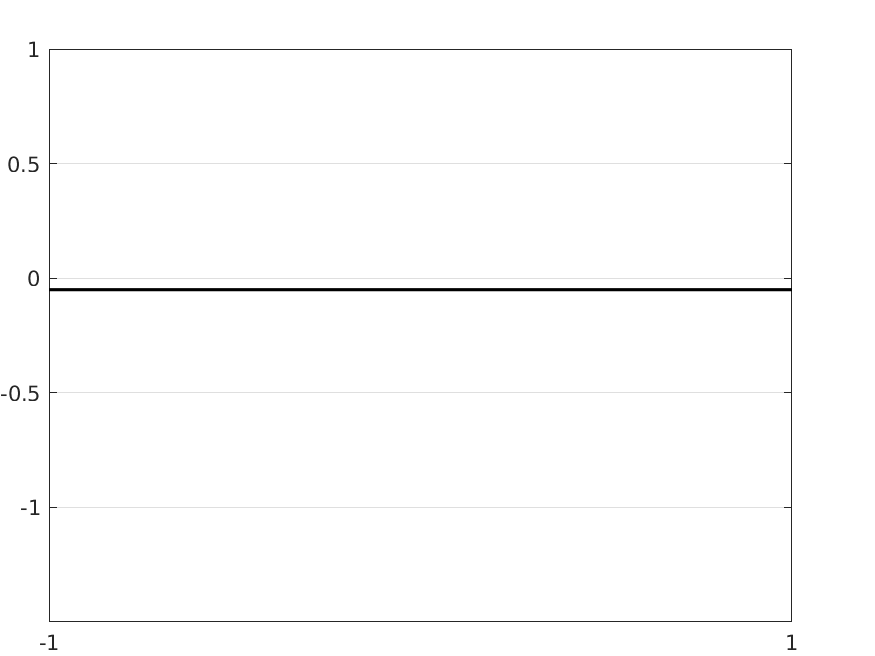}
                       & \includegraphics[width=.33\textwidth]{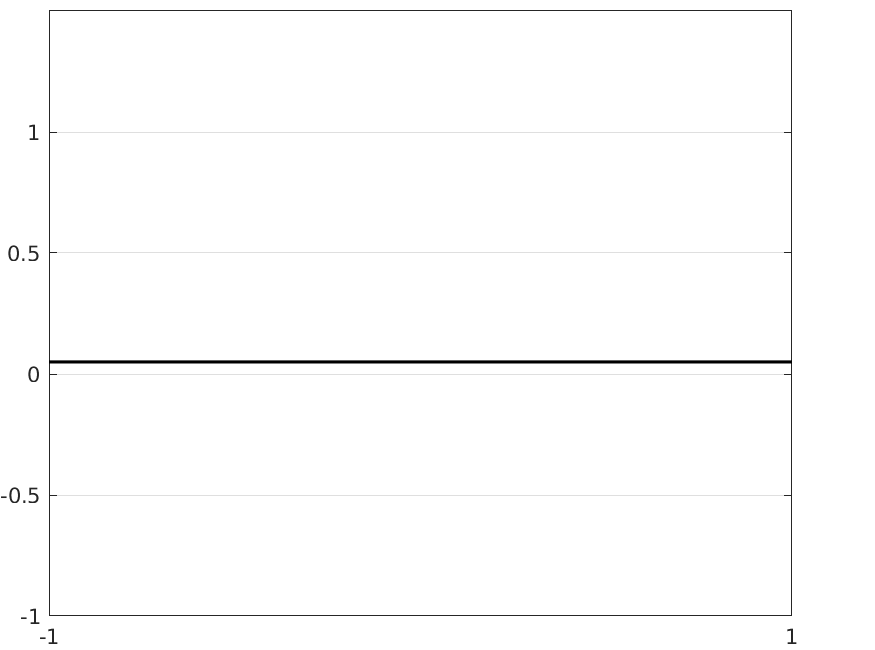}
                       & \includegraphics[width=.33\textwidth]{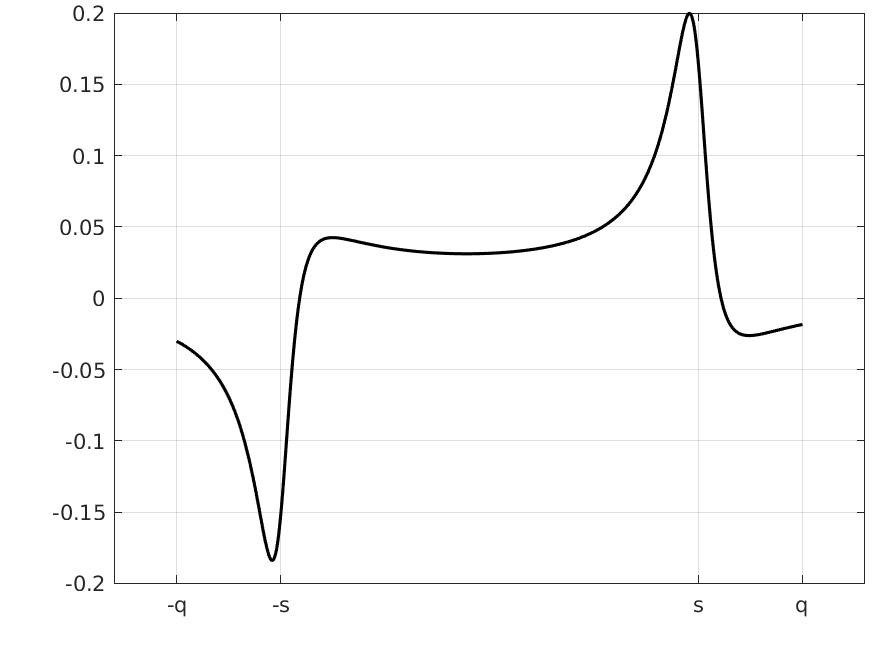}
                     \end{tabular}
    \caption{$m_1$, $m_2$ constant, associated $b_2[{\bm}]$.} 
    \label{figmcte}
        \end{figure}

\noindent 
Take $m_1 =-0.05$ on $S=[-1, 1]$ and $m_2 = 0.05$ on $S$ (whence ${\bm} = 0.05 \, \left({\be}_1-{\be}_2\right)$ and $m_i = 0.05 \, (-1)^i \, \chi_S$).  See Figure \ref{figmcte} and Table \ref{table3-magcte} for plots of  $m_i$ and of $b_2[{\bm}]$, computed as explained in Section \ref{sec:fp}. 

\begin{table}[htp]
    \centering
\begin{tabular}{|c|c|c|c|c|c|}
\hline
$\phi_i$ & $\lambda_i$ & $\langle m_1^{e} \rangle$& $\langle m_2^{e} \rangle$&$\varepsilon_1$ & $\varepsilon_2$\\[0.07cm]
\hline
in $L^2(K)$ & $10^{-5}$&$-0.1044$&$0.09581$&$4.4\cdot 10^{-4}$&$ 4.2\cdot 10^{-3}$\\ \hline
in $W^{1,2}_0(K)$& $10^{-8}$&$-0.0996$ &$0.0994 $& $3.8 \cdot 10^{-3}$ & $6.4\cdot 10^{-3}$\\ \hline
\end{tabular}

\caption{$\bm$ identically constant.}
  \label{table3-magcte}
\end{table}

\subsubsection{Magnetizations with large support}
\label{subs:mwls}

\begin{figure}[h]
      \begin{tabular}{lll}    \includegraphics[width=.3\textwidth]{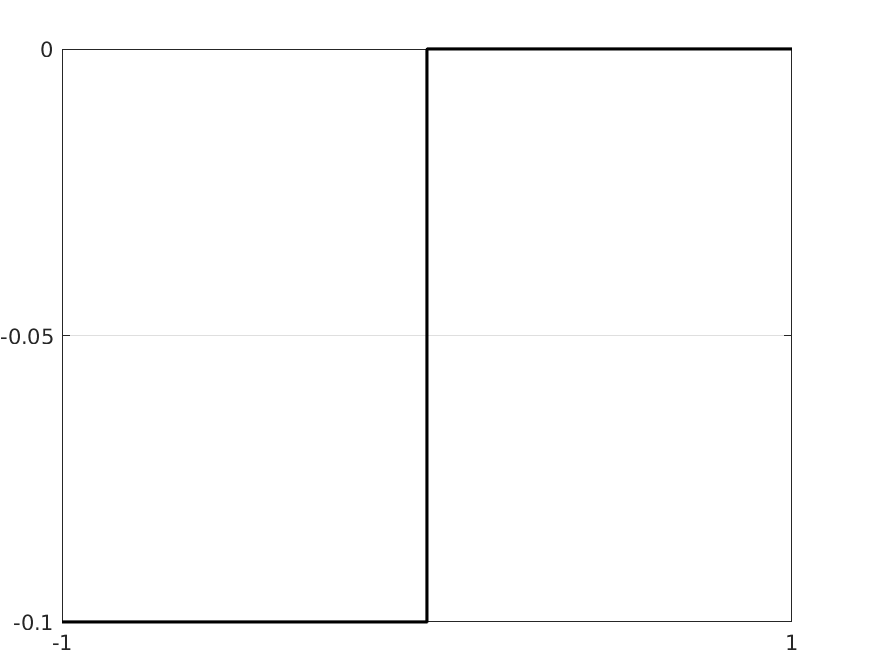}
        & \includegraphics[width=.33\textwidth]{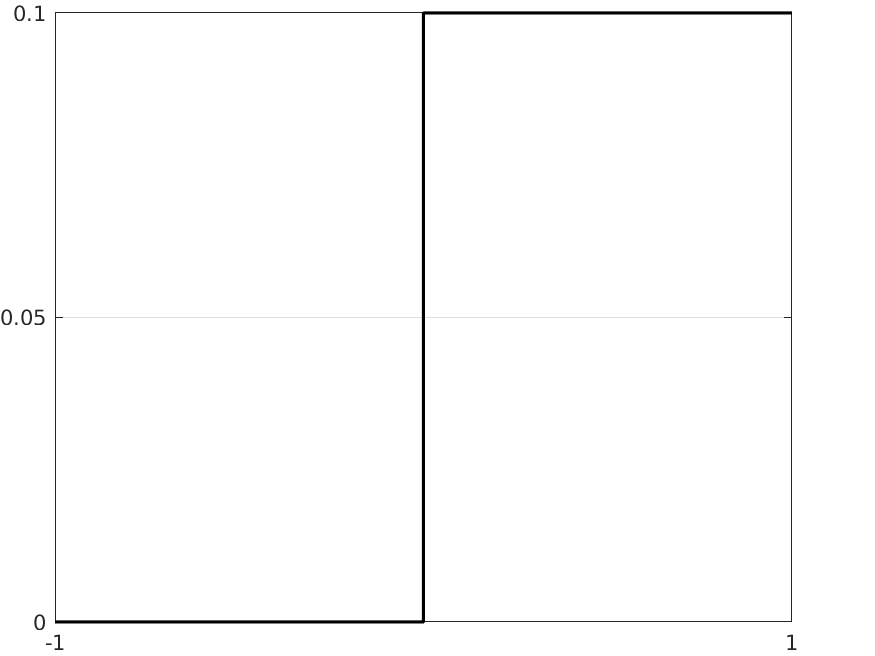}
        & \includegraphics[width=.33\textwidth]{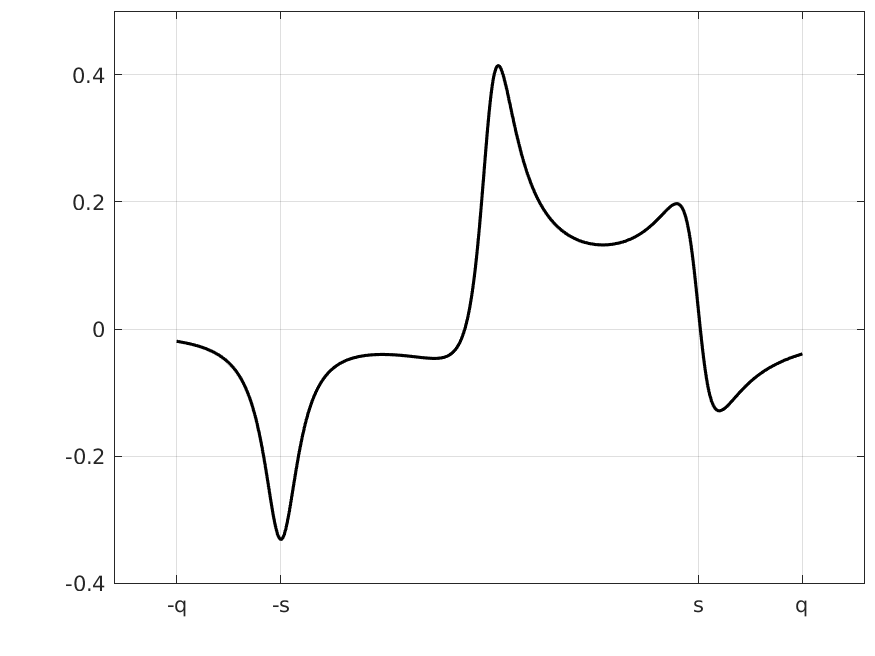}
            \end{tabular}
    \caption{$m_1$, $m_2$ with large support, associated $b_2[{\bm}]$.} 
    \label{figm0}
\end{figure}

\noindent 
Take $m_1 =-0.1$ on $[-1, 0]$, $0$ elsewhere in $S$, and $m_2 = 0.1$ on $[0, 1]$, see Figure \ref{figm0}, together with Table \ref{table2.2}.

\begin{table}[htp]
    \centering
\begin{tabular}{|c|c|c|c|c|c|}
\hline
$\phi_i$ & $\lambda_i$ & $\langle m_1^{e} \rangle$& $\langle m_2^{e} \rangle$&$\varepsilon_1$ & $\varepsilon_2$\\[0.07cm]
\hline
in $L^2(K)$ & $10^{-5}$&$-0.0999$&$0.0994$&$6.4\cdot 10^{-4}$&$ 5.5\cdot 10^{-3}$ \\ \hline
in $W^{1,2}_0(K)$& $10^{-8}$&$-0.1000$ &$0.0995 $& $4.4 \cdot 10^{-4}$ & $4.6\cdot 10^{-3}$\\ \hline
\end{tabular}

\caption{$\bm$ with large support.}
  \label{table2.2}
\end{table}

\subsubsection{Other magnetization}
\label{subs:msteps}

\begin{figure}[h]
  \begin{tabular}{lll}    \includegraphics[width=.3\textwidth]{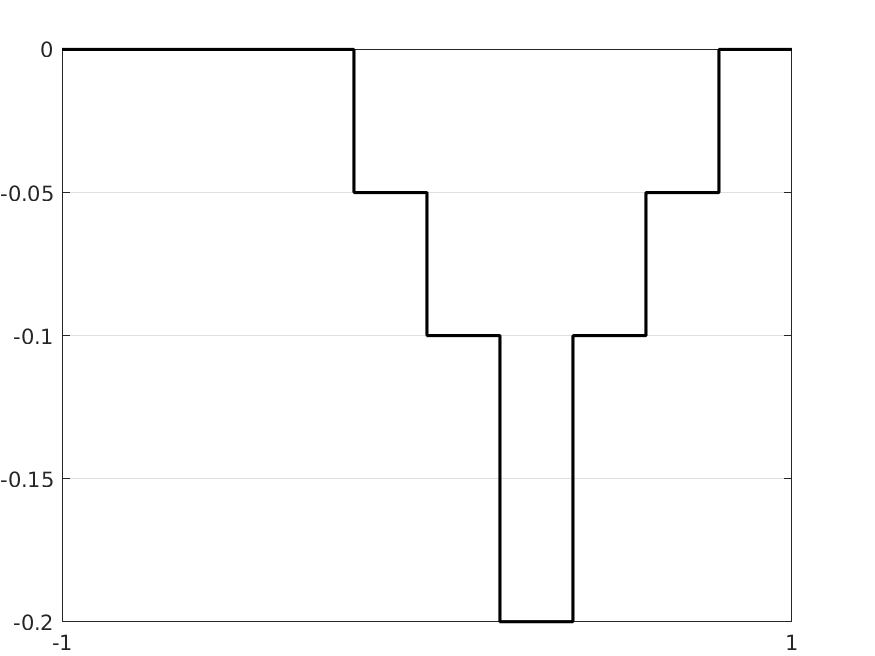} 
    & \includegraphics[width=.33\textwidth]{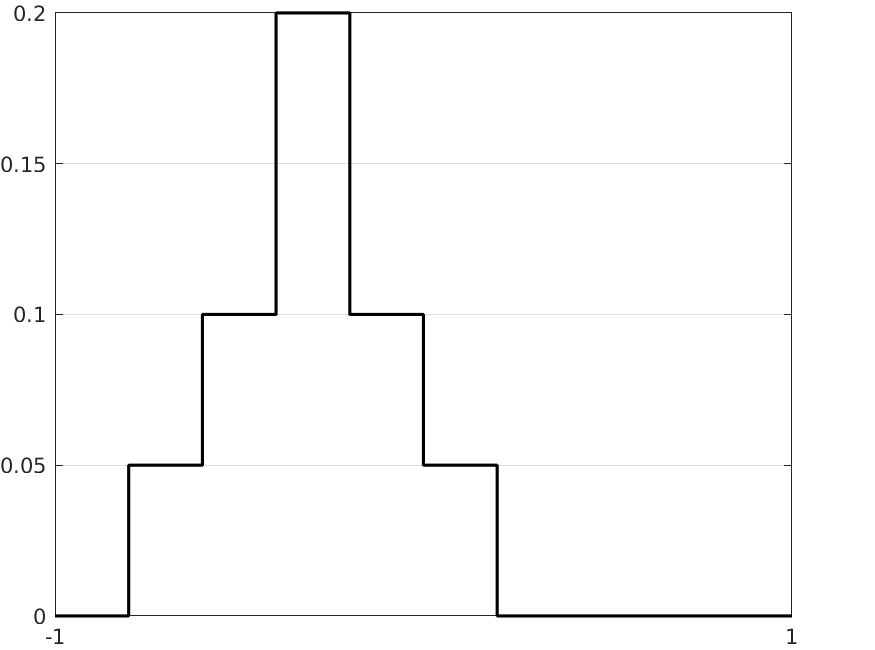}
    & \includegraphics[width=.33\textwidth]{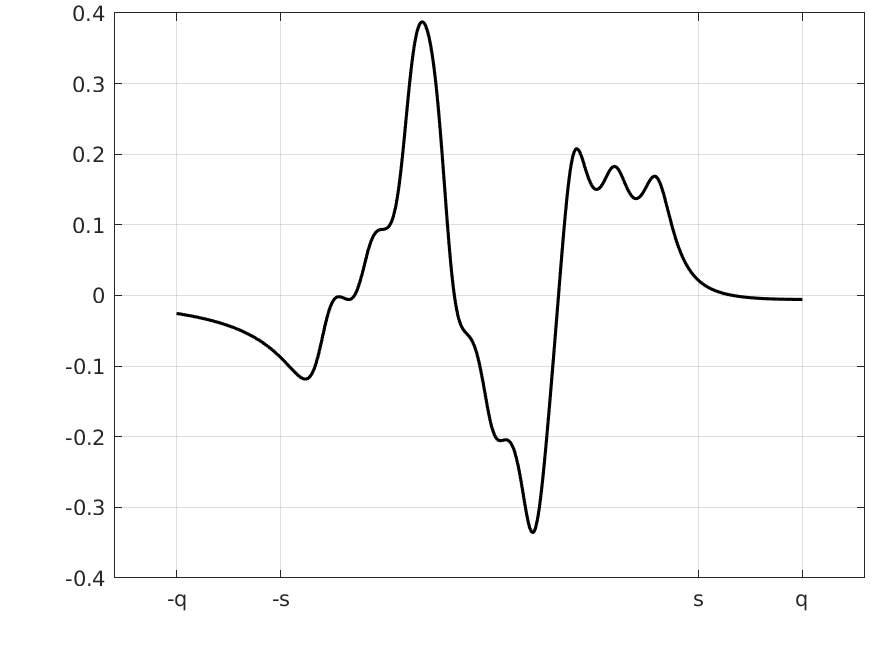}    \end{tabular}
    \caption{Others $m_1$, $m_2$, associated $b_2[{\bm}]$.} 
    \label{figmsteps}
\end{figure}

\noindent Here we take $m_1$ and $m_2$ as below. See Figure \ref{figmsteps} and Table \ref{table4-magsteps}.

$$\begin{matrix}
m_1=\begin{cases}
-0.05&\hbox{on } [-0.2, 0],\\
-0.1&\hbox{on } [0, 0.2],\\
-0.2&\hbox{on } [0.2, 0.4],\\
-0.1& \hbox{on }[0.4, 0.6],\\
-0.05&\hbox{on } [0.6, 0.8],\\
0&\hbox{elsewhere in }S,
\end{cases}&\qquad \qquad m_2=\begin{cases}
0.05&\hbox{on }[-0.8, -0.6],\\
0.1&\hbox{on } [-0.6, -0.4],\\
0.2&\hbox{on } [-0.4, -0.2],\\
0.1& \hbox{on }[-0.2, 0]\\
0.05&\hbox{on } [0, 0.2],\\
0&\hbox{elsewhere in }S.
\end{cases}
\end{matrix}
 $$

\begin{table}[htp]
    \centering
\begin{tabular}{|c|c|c|c|c|c|}
\hline
$\phi_i$ & $\lambda_i$ & $\langle m_1^{e} \rangle$& $\langle m_2^{e} \rangle$&$\varepsilon_1$ & $\varepsilon_2$\\[0.07cm]
\hline
in $L^2(K)$ & $10^{-5}$&$-0.0981$&$0.09855$&$1.9\cdot 10^{-2}$&$ 1.4\cdot 10^{-2}$ \\ \hline
in $W^{1,2}_0(K)$& $10^{-8}$&$-0.0977$ &$0.0989 $& $2.3 \cdot 10^{-2}$ & $1.1\cdot 10^{-2}$\\ \hline
\end{tabular}

\caption{Other $\bm$.}
  \label{table4-magsteps}
\end{table}

\subsubsection{Magnetizations with small support}
\label{subs:mwss}
\begin{figure}[h]
\begin{tabular}{lll} 
  \includegraphics[width=.33\textwidth]{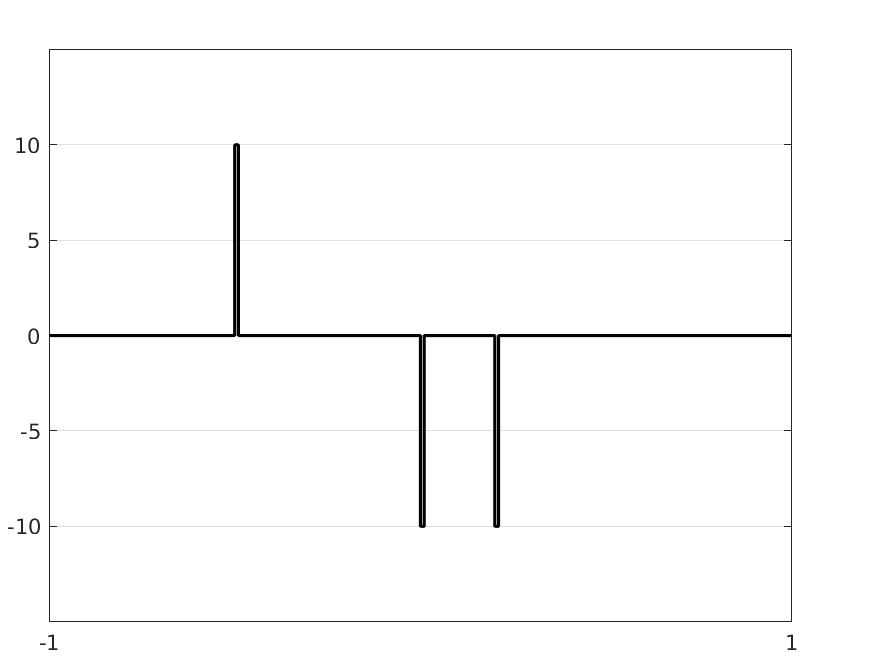}                  
                  & 
  \includegraphics[width=.33\textwidth]{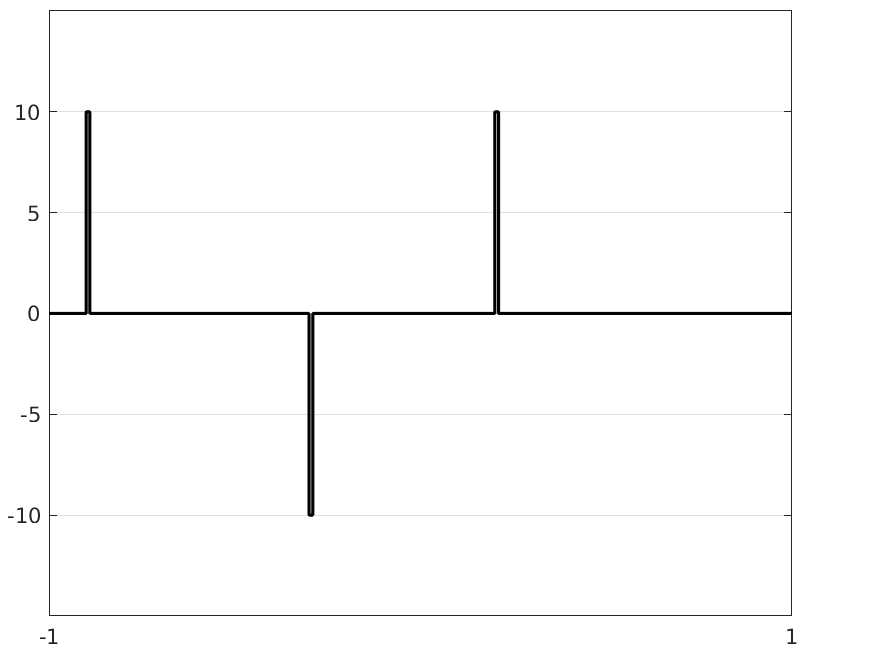}   
  &
\includegraphics[width=.33\textwidth]{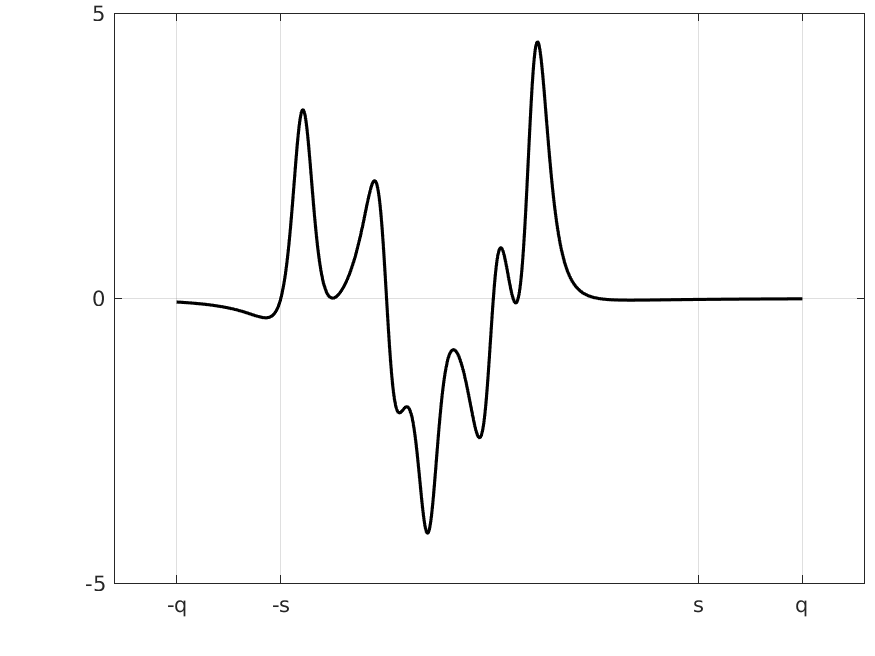}   
\end{tabular}
  \caption{$m_1$, $m_2$ with small support, associated $b_2[{\bm}]$.}  
  \label{figm}
\end{figure}

\noindent Here we take $m_1$ and $m_2$ as described below. See Figure \ref{figm} and Table \ref{table1.2bis}.

$$\begin{matrix}
m_1=\begin{cases}
10&\hbox{on } [-0.5, -0.49],\\
-10&\hbox{on } [0, 0.0.1],\\
-10&\hbox{on } [0.2, 0.21],\\
0&\hbox{elsewhere in }S,
\end{cases}&\qquad \qquad m_2=\begin{cases}
10&\hbox{on } [-0.9, -0.89],\\
-10&\hbox{on } [-0.3,-0.29],\\
10&\hbox{on } [0.2, 0.21],\\
0&\hbox{elsewhere in }S.
\end{cases}
\end{matrix}
 $$

      \begin{center} 
  \begin{table}[htp]
    \centering
\begin{tabular}{|c|c|c|c|c|c|} 
\hline
$\phi_i$ & $\lambda_i$ & $\langle m_1^{e} \rangle$& $\langle m_2^{e} \rangle$&$\varepsilon_1$ & $\varepsilon_2$\\[0.07cm] 
\hline

in $L^2(K)$ & $10^{-5}$&$-0.104$&$0.0958$& $ 4.4\cdot 10^{-2}$ & $4.2\cdot 10^{-2}$\\ \hline
in $W^{1,2}_0(K)$& $10^{-8}$&$ -0.1015$&$0.0969$& $1.5 \cdot 10^{-2}$ & $3.1\cdot 10^{-2}$\\ \hline 
\end{tabular}

\caption{$\bm$ with small support.}
\label{table1.2bis}
  \end{table}
      \end{center}

\subsubsection{Comments, discussion}
\noindent
Overall, as presented in the above tables, we obtain quite accurate results for  net moment estimation. The net moments of magnetizations with large support are more precisely estimated than those of magnetizations with small support: the value $\varepsilon_i$ of the error decreases whenever the size of the support increases, for a same value of the parameter $\lambda_i$. This phenomenon is mostly true in $L^2(K)$ for the present examples. These examples may not always provide smaller estimation errors $\varepsilon_i$ in $W^{1,2}_0(K)$ than in $L^2(K)$, since the behaviour of $\varepsilon_i$ depend on the specific data (smoothness properties of $\bm$, size of its support, ...) and parameters. These properties, together with the influence of the noise in the computations and in the data, will be studied  in a forthcoming work.

\section{Perspectives, conclusions}
\label{sec:pc}

\noindent In order to complete the results of Section \ref{sec:mo}, some properties of $b_2$ remain to be studied. Even if $b_2$ is injective, it is not coercive (strongly injective), in the sense that  $\left\|b_2[{\bm}]\right\|_{L^2(K)}$ can  be small even when $\left\|{\bm}\right\|_{L^2(S, \RR^2)}$ is not small: there may exist such ``almost silent'' source terms ${\bm}$. 
However, lower bounds for $\left\|b_2\right\|$ can be established if we restrict to truncated Fourier expansions, showing that these cannot be ``almost silent''.
Also, one can estimate the constants involved in the upper bounds of $\left\|b_2\right\|$ and $\left\|b_2^*\right\|$.
The  spectral study of the operator $b_2 \, b_2^*$ remains to be fulfilled, in view of the numerical analysis of the BEP.\\

 \noindent The BEP studied in Section \ref{sec:beps} can be stated with slightly more general constraints.
In particular,  let $f\in L^2(K)$. Consider the following BEP in $L^2(K)$:

\[
\min_{\phi\in L^{2}(K)\atop{
\Vert \phi-f\Vert_{L^{2}(K)}\leq M}}\Vert b^*_{2}[{\phi}] -{\be}\Vert_{L^2(S,\RR^2)} \, \label{BEPf}\tag{BEP$_f$}.
\]

\noindent One can prove 
existence and uniqueness of the function $\phi_{o}\in L^2(K)$ solution of \eqref{BEPf} using the same arguments as in the proof of Proposition \ref{prop:wpBEP}. Moreover, if
 $ {\be}  \not \in  \text{Ran}\,b_2^*$ (or if ${\be}= b_2^{*}[\phi]$ for some  $\phi \in L^2(K)$ such that $\Vert \phi - f\Vert_{L^2(K)}\geq M$), then the solution $\phi_{o}$ still saturates the constraint, i.e. $\Vert\phi_{o}-f\Vert_{L^2(K)}=M$, and is given by the following implicit equation, with $\gamma  < 0$, see \eqref{BEP}:

$$b_2b_2^*[\phi_{o}]-\gamma \phi_{o}=b_2[{\be}]+\gamma f \, .$$

\noindent Let now $f\in W_0^{1,2}(K)$. The BEP can be stated as follows:
\[
\min_{\phi\in W_0^{1,2}(K)\atop{
\Vert \phi-f\Vert_{W_0^{1,2}(K)}\leq M}}\Vert b^*_{2}[{\phi}] -{\be}\Vert_{L^2(S,\RR^2)} \, \label{BEPfW}\tag{BEP$_{f,W}$}.
\]

\noindent Again, as for \eqref{BEP1}, a solution $\phi_o\in W_0^{1,2}(K)$ exists, is unique, and saturates the constraint if
$ {\be}  \not \in  b^*_{2|_W}[W_0^{1,2}(K)]$ (or if ${\be} =b_{2|_W}^*[\phi]$ for $\phi \in {L^2(K)}$ such that $\Vert \phi - f\Vert_{W_0^{1,2}(K)}\geq M$):~$\Vert \phi_o-f \Vert_{W_0^{1,2}(K)}= M$. It satisfies the following equation, with $\gamma < 0$:

$$ b_2b_2^*[\phi_{o}]+\gamma \phi''_{o}=b_2[{\be}]+\gamma f \, .$$

\noindent The above problem will be considered for $\phi_{o} \in W^{1,2}(K)$ with norm constraint therein, where it is also well-posed and gives promising numerical results. Moreover, $W^{1,2}(K)$ functions are continuous on $K$ in the present 1-dimensional setting, which is a suitable property of the solution, in view of computing its inner product with the available measurements in a pointwise way. Note however that the continuity of the solution is already granted by the critical point equation \eqref{CPEmo} above and the smoothness properties of functions in $\text{Ran} \, b_2$. 

\noindent Some BEP of mixed type, with  constraints in other norms, like uniform, $L^1$, BMO, could also be formulated, as well as related extremal problems consist in looking for the best bounded extension of a given data $b_d \simeq b_2[{\bm}]$ from $K$ to $\RR$.\\

\noindent  Following Section \ref{sec:comput}, further numerical analysis and improvements of the computational schemes will be considered.
In particular, we expect that refined implementations of the computations in the Fourier domain will provide more efficient and accurate recovery schemes, in the present setting as well as in the 3D one, see \cite{BEP-moments}. Both a quantitative and a qualitative study of the relations $\lambda \mapsto M(\lambda)$ and $\lambda \mapsto \epsilon_i(\lambda)$ (for the different error terms) remain to be done. 
This will be the topic of a further work, 
together with a study concerning the influence of the parameters $s$, $q$, $h$, and of the characteristics of $m_i$, on the behaviour of the solutions. \\

\noindent Both explicit and asymptotic expressions, in terms of the size $q$ of the measurement set $K$, relating the first moments of ${\bm}$ to those of $b_2[{\bm}]$, could also be derived, as in \cite{asympt} for the 3D situation.  Besides, the quantities $b^*_{2}[x^k]$ for monomials $x^k$ could be exactly computed. %

\noindent Other functions $\be$ may be considered as well, both for the present net moment estimation problem or for higher order or local moments recovery. One can also take the components $e_i$ to be appropriate band-limited basis functions, of 
Slepian type  \cite{Ledford}.

\noindent Particularly interesting magnetizations are the unidirectional ones, of which the moment recovery using the above linear estimator process still requires a specific study. Also, more general situations where  ${\bm}$ belongs to $ L^1(S, \RR^2)$ or is a measure remain to be considered, see \cite{Betal-IP}, together with magnetizations supported on a two-dimensional set.\\

\noindent  Finally, the full inversion problem of recovering ${\bm}$ in $L^2(S, \RR^2)$ itself from $b_2[{\bm}]$ will be considered in a further work. It can also be stated as a BEP of  which the solution will involve the operator $b_2 \, b_2^*$. Again more general magnetization distributions will be studied.

\section*{Acknowledgments}
\noindent The authors warmly thank Jean-Paul Marmorat for his help in numerical implementations.

\newpage

\appendix
\section*{Appendix: numerical illustrations}

\vfill
\end{document}